\documentclass[11pt, a4paper]{article}
\usepackage{graphicx}
\usepackage{color}
\usepackage{amsmath}
\usepackage{amssymb}
\usepackage{amscd}
\usepackage{bbm}
\usepackage{setspace}
\usepackage{slashbox}
\usepackage{tikz}

\usepackage{marginnote}

\usepackage{hyperref}
\hypersetup{
    bookmarks=true,         
    colorlinks=true,       
    linkcolor=red,          
    citecolor=green,        
    filecolor=magenta,      
    urlcolor=cyan           
}

\newcommand{\R}{\mathbb{R}}
\newcommand{\inr}[1]{\bigl< #1 \bigr>}

\newcommand{\E}{\mathbb{E}}

\newcommand{\eps}{\varepsilon}

\newcommand{\F}{{\cal F}}

\newcommand{\cA}{{\mathcal{A}}}
\newcommand{\bX}{{\mathbb{X}}}

\newcommand{\hf}{{\hat f_N}}

\newcommand{\cB}{{\mathcal{B}}}
\newcommand{\cH}{{\mathcal{H}}}
\newcommand{\cD}{{\mathcal{D}}}
\newcommand{\cX}{{\mathcal{X}}}
\newcommand{\cF}{{\mathcal{F}}}
\newcommand{\cL}{{\mathcal{L}}}
\newcommand{\norm}[1]{\left\|#1\right\|}%
\newcommand{\cN}{{\mathcal{N}}}

\newcommand{\Pro}{\mathbb{P}}

\DeclareMathOperator*{\argmin}{argmin}

\newtheorem{Theorem}{Theorem}[section]
\newtheorem{Lemma}[Theorem]{Lemma}

\newtheorem{Definition}[Theorem]{Definition}

\newtheorem{Remark}[Theorem]{Remark}

\numberwithin{equation}{section}

\def \proof {\noindent {\bf Proof.}\ \ }

\def \endproof
{{\mbox{}\nolinebreak\hfill\rule{2mm}{2mm}\par\medbreak}}
\def\IND{\mathbbm{1}}

\begin{document}
\title{{Learning subgaussian classes : Upper and minimax bounds}}
\author{Guillaume Lecu\'e${}^{1,3, 5}$  \and Shahar Mendelson${}^{2,4,6}$}

\footnotetext[1]{CNRS, CREST, ENSAE, 3 avenue Pierre Larousse, 92240 Malakoff, France}
\footnotetext[2]{Department of Mathematics, Technion, I.I.T, Haifa
32000, Israel.}
\footnotetext[3] {Email:guillaume.lecue@ensae.fr }
\footnotetext[4] {Email:shahar@tx.technion.ac.il}
\footnotetext[5]{Supported by French National Research Agency (ANR) under the grants Labex Ecodec (ANR-11-LABEX-0047) and by the "Chaire Economie et Gestion des Nouvelles Donn\'ees", under the auspices of Institut Louis Bachelier, Havas-Media and Paris-Dauphine.}
\footnotetext[6]{Supported by the Mathematical Sciences Institute, The Australian National University and by the Israel Science Foundation.}

\maketitle

\section*{Preface}
Most the results contained in this note have been presented at the SMF meeting, which took place in May 2011; the rest have been obtained shortly after the time of the meeting.

The question we study has to do with the optimality of Empirical Risk Minimization as a learning procedure in a convex class -- when the problem is subgaussian. Subgaussian learning problems are a natural object because they are the simplest unbounded learning scenarios. However, an additional reason for studying such problems was that at the time of the SMF meeting, the technical machinery required for the analysis of more heavy-tailed problems was simply not known. Since 2011, significant progress has been made in the understanding of learning problems in heavy-tailed situations \cite{MR3367000,LWCG,LM_sparsity,LugMen}, though this progress does not make the results presented here obsolete. We show that ERM performed in a convex class is an optimal learning procedure (in a sense that will be clarified) when the learning problem is subgaussian. This happens to be a rather special feature of subgaussian learning problems, and under weaker tail assumptions ERM fails to deliver the optimal accuracy/confidence trade-off at the high level of accuracy we are interested in here.

The results presented here are complemented in \cite{MenLvG}, which also focuses on subgaussian learning problems and addresses some of the cases that have not been resolved in this note.

\section{Introduction and main results} \label{sec:intro}
Let $\cD:=\{(X_i,Y_i):i=1,\cdots,N\}$ be a set of $N$ i.i.d random variables with values in $\cX\times \R$. From a statistical standpoint, each $X_i$ can be viewed as an input associated with a real-valued output $Y_i$.  Given a new input $X$, one would like to guess its associated output $Y$, assuming that $(X,Y)$ is distributed according to the same probability distribution that generates the data $\cD$. To that end, one may use $\cD$ to construct a function $\hat f_N(\cD,\cdot)=\hat f_N(\cdot)$,  and the hope is that  $\hat f_N(X)$ is close to $Y$ in some sense.

Here, we will consider the \textit{squared loss function} $\ell:\R\times \R\mapsto\R$, defined by $\ell(u,v)=(u-v)^2$, as a way of measuring the pointwise error $\ell(f(X),Y)$, and the resulting \textit{squared risk} is
\begin{equation*}
  R(f)=\E \big(f(X)-Y\big)^2 \mbox{ and } R(\hat f_N)=\E\big(\big(\hat f_N(X)-Y\big)^2|\cD\big).
\end{equation*}

In the classical statistics setup, one usually assumes that the regression function of $Y$ given $X$ belongs to some particular function space (called a \textit{statistical model}). In contrast, in the learning setup on which we focus here, one is given a function class $\cF$ (sometimes, called a model as well), and the goal is to construct a procedure $\hat f_N$ that satisfies a \textit{sharp} or \textit{exact oracle inequality} (following \cite{TsyCOLT07}; such bounds are called excess risk bounds in \cite{MR2319879} and \cite{MR2829871}). An exact oracle inequality ensures that with high probability,
\begin{equation}\label{eq:Exact-Oracle-Ineq}
  R(\hat f_N)\leq \inf_{f\in \cF}R(f)+{\rm residue},
\end{equation}
and one would like to make the residue in \eqref{eq:Exact-Oracle-Ineq} as small as possible.

For the sake of simplicity, we assume that there is some $f^*\in\cF$ minimizing the risk in $\cF$ (though the claims presented here remain true even without that assumption), and we set
\begin{equation*}
f^*\in\argmin_{f\in\cF}R(f).
\end{equation*}
Note that in \eqref{eq:Exact-Oracle-Ineq} the performance of the procedure $\hf$ is compared to the best performance possible in $\cF$, i.e., to the risk of the best element $f^* \in \cF$. This exhibits the point of view of Learning Theory, where one wishes to identify a function that is almost as good as the best possible in $\cF$, regardless of whether the best function in $\cF$ has a small risk. It is different from typical questions in classical Statistics, where a statistical model is given and the risk of an estimator is compared to the one of the regression function (or Bayes rule). The latter are usually called \textit{excess risk bounds} (cf. \cite{massart03}) and are actually very different from exact oracle inequalities like \eqref{eq:Exact-Oracle-Ineq} (see, for example, \cite{MR2933668} or Chapter~1.3 in \cite{HDR-lecue} for more details on those differences).

The performance of a procedure is measured relative to a set of admissible targets $Y$ in  some class of random variables ${\cal Y}$. Naturally, one would like to make ${\cal Y}$ as large as possible, for example, all random variables $Y$ bounded by $1$, all the random variables $Y$ in $L_p$ for some $p >2$, or a similar weak condition of that flavor.

\begin{Definition}\label{def:accuracy_deviation}
Let $\hat{f}_N$ be a learning procedure, that is, a map from the set $(\Omega \times \R)^N$ into $\cF$. Let $0<\delta_N<1$ and $\eps_N>0$. We say that $\hf$ performs with \textbf{accuracy} $\eps_N$ and \textbf{confidence} $1-\delta_N$ relative to the set of admissible targets ${\cal Y}$ if for any $Y \in {\cal Y}$, $R(\hat f_N)\leq \inf_{f\in \cF}R(f)+\eps_N$ with probability larger than $1-\delta_N$, and the probability is measured with respect to the product measure endowed by the joint distribution of $X$ and $Y$.
\end{Definition}

Clearly, while the true risk of $f$ is not known, simply because $X$ and $Y$ are not known, one still has access to its empirical counterpart:
\begin{equation*}
  R_N(f)=\frac{1}{N}\sum_{i=1}^N \big(f(X_i)-Y_i\big)^2.
\end{equation*}
Thus, a natural procedure that comes to mind is finding a function in $\F$ that best fits the data: a minimizer of the empirical risk in $\cF$. This procedure is called \textit{empirical risk minimization} (\textit{ERM}) and is defined by
\begin{equation*}
  \hat f\in\argmin_{f\in \cF}R_N(f).
\end{equation*}
ERM has been studied extensively over the last $40$ years (see, e.g. \cite{vapnik98}, \cite{MR2319879}, \cite{MR2829871} and references therein), and the main goal has always been to identify connections between the structure of $\cF$ and the accuracy and confidence that ERM yields, while trying to minimize the restrictions on ${\cal Y}$. Among the natural questions regarding the performance of ERM are:
\begin{description}
\item{1.} Given any confidence parameter $0<\delta_N<1/2$, what is the error rate $\eps_N$ that one may obtain using ERM, and what features of $\cF$ govern that rate?
\item{2.} Given any $0<\delta_N<1/2$, is ERM an optimal procedure for the confidence level $\delta_N$? In other words, is there a procedure that can perform with a better accuracy than ERM, given the same confidence level?
\end{description}

The majority of results on the performance of ERM have been obtained in the bounded case: when $\sup_{f \in \cF} |\ell(Y,f(X))|\leq b$ almost surely, or, alternatively, when the envelope function $\sup_{f\in\cF}|\ell(Y,f(X))|$ is well behaved in some weaker sense (e.g., has a sub-exponential tail). A result in this direction is from \cite{MR2166554} (see Corollary~5.3 there) which we formulate using the notation of Theorem~5.1 in \cite{MR2829871}.

For any $\gamma>0$, let
\begin{equation}\label{eq:fixed_point_rad_comp}
k_N(r)=\E\sup\Big(\Big|\frac{1}{\sqrt{N}}\sum_{i=1}^N\eps_i (f-f^*)(X_i)\Big|:f\in\cF,\norm{f-f^*}_{L_2(\mu)}\leq 2r\Big),
\end{equation}
and set $k_N^*(\gamma)=\inf\left\{r>0:8k_N(r)\leq \gamma r^2\sqrt{N}\right\}$.
\begin{Theorem} \label{thm:Kolt}
There exist absolute constants $c_0,c_1$ and $q>2$ for which the following holds. If ${\cal Y}$ consists of functions that are bounded by $1$ and $\cF$ is a convex class of functions that are bounded by $1$, then for any $Y \in {\cal Y}$ and every $t>0$, with probability at least $1-c_0\exp(-t)$,
\begin{equation}
  \label{eq:koltchinskii}
  R(\hat f)\leq \inf_{f\in \cF}R(f)+c_1\max\Big\{\big(k_N^*(1/q)\big)^2,\frac{t}{N}\Big\}.
\end{equation}
\end{Theorem}

A result of a similar flavor was obtained in \cite{MR1240719}: let $N(A,B)$ be the number of translates of $B$ needed to cover $A$. Set $D$ to be the unit ball in $L_2(\mu)$ and let
\begin{equation}
  \label{eq:Birge-massart-93}
  \sigma^*=\inf\left\{r>0:\int_{c_0 r^2}^{c_1 r}\log^{1/2}N\big(\cF\cap(2rD),\eps D\big)d\eps\leq c_2 r^2\sqrt{N}\right\},
\end{equation}
for absolute constants $c_0,c_1,c_2$.

The result in \cite{MR1240719} is that under various assumptions on the class $\cF$ (assumptions that allow one to upper bound the function $k_N(r)$ using the entropy integral in \eqref{eq:Birge-massart-93}), $(\sigma^*)^2$ may serve as a residual term.

These two facts rely heavily on the assumption that $\cF$ and ${\cal Y}$ are bounded in $L_\infty$ and their proofs do not extend beyond the bounded case.

Our aim here is to study unbounded problems and without any assumption on the envelope of $\{\ell(f(X),Y):f\in\cF\}$. The next natural step is the subgaussian framework, as it captures many typical applications in which the functions involved are unbounded: for example, regression with a gaussian noise; compressed sensing; matrix completion; phase recovery, etc. (see \cite{MR2807761,MR2723472,IEEE-CT,MR2382644,Gross,MR1093464}), all of which have been studied in the subgaussian framework.

\begin{Definition}\label{def:subgaussian-class}
Let $\mu$ be a probability measure and let $X$ be distributed according to $\mu$. The $\psi_2(\mu)$-norm of a function $f$ is
$$
\norm{f}_{\psi_2(\mu)}=\inf\left\{c>0:\E\exp(f^2(X)/c^2)\leq 2\right\}.
$$
The space of functions with a finite $\psi_2$-norm is denoted by $L_{\psi_2}=L_{\psi_2(\mu)}$.

A function class  $\cF \subset L_2(\mu)$ is \textbf{$L$-subgaussian} with respect to the probability measure $\mu$ if for every  $f,h \in \cF \cup \{0\}$, $\norm{f-h}_{\psi_2(\mu)}\leq L \norm{f-h}_{L_2(\mu)}$.
\end{Definition}

Note that for any $f\in L_{\psi_2}$, $\norm{f}_{L_2(\mu)}\leq\norm{f}_{\psi_2(\mu)}$. A class is a subgaussian class when the reverse inequality holds, and in particular when the $\psi_2$ and $L_2$ norms are equivalent on $\cF$. 

Note that norm equivalence is very different from being bounded. Having such a norm equivalence implies that $|f| \sim \|f\|_{L_2(\mu)}$ on a relatively large event. In contrast, even though a bounded function has a finite $\psi_2$ norm (by selecting $c \sim \|f\|_{L_\infty}$ in the definition of the $\psi_2$ norm), the fact that $f$ is bounded does not mean that $\|f\|_{\psi_2}$ is equivalent to $\|f\|_{L_2}$, nor that $|f| \sim \|f\|_{L_2(\mu)}$ on a relatively large event. Because of the substantial difference between the two notions, one should not expect that learning procedures exhibit the same performance when one assumes that $\cF$ is bounded in $L_\infty$ or when the $\psi_2(\mu)$ and $L_2(\mu)$ norms are equivalent on $\cF$.

\vskip0.3cm

Let us turn to some examples of subgaussian classes of functions. Probably the most interesting collection of examples that belong to the subgaussian framework is classes of linear functionals on $\R^d$.

\begin{Definition}
  \label{def:subgaussian-measure}
A probability  measure $\mu$ on $\R^d$ is \textbf{$L$-subgaussian}, if for every $t\in \R^d$, $\norm{\inr{t,\cdot}}_{\psi_2(\mu)}\leq L\norm{\inr{t,\cdot}}_{L_2(\mu)}$. The measure  $\mu$ is \textbf{isotropic} if $\norm{\inr{t,\cdot}}_{L_2(\mu)}=\norm{t}_{\ell_2^d}^2$ for every $t\in\R^d$, where $\norm{\cdot}_{\ell_2^d}$ denotes the Euclidean norm in $\R^d$.
\end{Definition}

There are many natural examples of subgaussian measures on $\R^d$:
\vskip0.3cm
\noindent{$\bullet$} Let $x$ be a real-valued random variable that has mean-zero and variance $1$. If $\|x\|_{\psi_2(\mu)} \leq L \|x\|_{L_2(\mu)}$ and $x_1,\ldots,x_d$ are independent copies of $x$, then it is straightforward to verify that for every $a \in \R^d$,
$$
\Big\|\sum_{i=1}^d a_ix_i\Big\|_{\psi_2(\mu)} \lesssim L\Big\|\sum_{i=1}^d a_ix_i\Big\|_{L_2(\mu)},
$$
where here, and throughout this note we write $u\lesssim v$ if $u\leq c_0 v$ for an absolute constant $c_0$. Thus, the measure associated with the random vector $X=(x_1,...,x_d)$ is $cL$-subgaussian. Also, the measure is clearly isotropic.

Natural examples of such product measures are the uniform measure on the combinatorial cube $\{-1,1\}^d$, the uniform measure on the cube $[-1,1]^d$ or the canonical gaussian measure in $\R^d$.
\vskip0.5cm
\noindent{$\bullet$}
Let $2 \leq p < \infty$ and denote by $B_p^d$ the unit ball of $(\R^d, \| \ \|_{\ell_p})$. The uniform probability measure on $d^{1/p}B_p^d$ is $L$-subgaussian for an absolute constant $L$ (see~\cite{MR2123199}), despite the fact that its coordinates are not independent.
\vskip0.5cm
\noindent{$\bullet$} Let $X=(x_i)_{i=1}^d$ be an \textit{unconditional} random vector (that is,  $(\eps_i x_i)_{i=1}^d$ has the same distribution as $X$ for every choice of signs $(\eps_i)_{i=1}^d$). If $\E x_i^2 \geq c^2$ for every $1 \leq i \leq d$ and $X$ is supported in $R B_\infty^d$, then it is $L$-subgaussian for $L \lesssim R/c$. Indeed, one may show that for every $f \in L_{\psi_2(\mu)}$,
$$
c_1\|f\|_{\psi_2(\mu)} \leq  \sup_{p \geq 2} \frac{\|f\|_{L_p(\mu)}}{\sqrt{p}} \leq c_2\|f\|_{\psi_2(\mu)}
$$
for suitable absolute constants $c_1$ and $c_2$ (see, for instance, Corollary~1.1.6 in \cite{MR3113826}). Thus, it suffices to verify that for every $t  \in \R^d$ and every $p \geq 2$,
$$
\|\inr{t,\cdot}\|_{L_p(\mu)} \leq L\sqrt{p}\|\inr{t,\cdot}\|_{L_2(\mu)}.
$$
By Khintchine's inequality (see, for example, \cite{LT:91}),
\begin{align*}
\|\inr{X,t}\|_{L_p}^p =   \E \Big|\sum_{j=1}^d x_j t_j\Big|^p = \E_X \E_\eps \Big|\sum_{j=1}^d \eps_j x_j t_j \Big|^p
\lesssim   p^{p/2} \E_X  \Big(\sum_{j=1}^d x_j^2 t_j^2\Big)^{p/2} \lesssim p^{p/2} R^p\|t\|_{\ell_2^d}^p.
\end{align*}
Also,
$$
\|\inr{X,t}\|_{L_2}^2 = \E_X \E_\eps \Big(\sum_{i=1}^d \eps_i x_i t_i\Big)^2 = \E_X \sum_{i=1}^d x_i^2t_i^2 \geq c^2\|t\|_{\ell_2^d}^2,
$$
proving the claim.
\vskip0.5cm
\noindent{$\bullet$} If $x$ is a mean-zero, variance one, $L$-subgaussian random variable, and  $X=(x_{i,j})$ is a matrix whose coordinates are independent copies of $x$, then $X$ defines a $cL$ subgaussian, isotropic  measure on the space of matrices of the right dimensions, relative to the natural trace inner product. The same holds if $X$ has independent rows, distributed according to an isotropic, $L$-subgaussian random vector. The proof of both facts is straightforward and are omitted.

\vskip0.3cm

These examples show that even the seemingly restricted setup of classes of linear functionals on $\R^d$ endowed with an $L$-subgaussian measure is encountered in many natural (and well studied) examples.

\vskip0.3cm

 The strategy we use here for the study of ERM is the isomorphic method, introduced in \cite{MR2240689} and analyzed there in the bounded setup. Before presenting it, recall that the excess loss of $f$ is
\begin{equation}
  \label{eq:Excess-loss}
  \cL_f(x,y)=\ell(f(x),y)-\ell(f^*(x),y)=(f(x)-y)^2-(f^*(x)-y)^2
\end{equation}
and set
$$
P\cL_f=\E \cL_f(X,Y) \ \ {\rm and} \ \ P_N \cL_f=\frac{1}{N}\sum_{i=1}^N \cL_f(X_i,Y_i).
$$
A rather obvious but very useful observation is that for every $f \in \cF$, $P \cL_f \geq 0$, while the empirical minimizer $\hat{f}$ satisfies that $P_N \cL_{\hat{f}} \leq 0$.

\vskip0.3cm

The isomorphic method is based on the following idea. Consider an event $\Omega_0$, on which for every function $f$ in the set $\{f\in\cF:P\cL_f\geq \lambda_N\}$,
\begin{equation}\label{eq:iso-property}
  \frac{1}{2}P \cL_f\leq P_N\cL_f\leq \frac{3}{2}P\cL_f.
\end{equation}
It follows that on $\Omega_0$, ERM produces $\hat{f}$ that satisfies
$$
R(\hat f)\leq \inf_{f\in\cF}R(f)+\lambda_N,
$$
because $P_N\cL_{\hat{f}} \leq 0$; therefore, $\hat{f} \not \in \{f \in \cF : P{\cL_f} \geq \lambda_N\}$.

Consequently, an exact oracle inequality with a confidence parameter $\delta_N$ may be derived by identifying $\lambda_N$ for which $\Omega_0$ has probability at least $1-\delta_N$; that is, the level $\lambda_N$ for which
\begin{equation*}
  \sup_{\{f\in\cF:P\cL_f\geq \lambda_N\}}\Big|\frac{1}{N}\sum_{i=1}^N\frac{\cL_f(X_i,Y_i)}{P\cL_f}-1\Big|\leq \frac{1}{2}
\end{equation*}
with probability at least $1-\delta_N$ (see Theorem~4.4 in \cite{MR2829871} for results of a similar flavor).

\begin{Remark}
Note that only the lower estimate in \eqref{eq:iso-property} is needed for the argument outlined above to work. This observation is the key in the application of the recent works on the {\it small-ball method} in learning theory (cf. \cite{MR3367000}), which allows one to deal with heavy-tailed scenarios that are far more general than subgaussian problems.
\end{Remark}

Just like $k_N^*$ in \eqref{eq:koltchinskii} and $\sigma^*$ in \eqref{eq:Birge-massart-93} -- and many other well known estimates on the performance of ERM (e.g. \cite{van_de_geer00,MR2829871,MR2319879}) -- the residual term we use is defined in terms of fixed points. Unlike $k_N^*$ and $\sigma^*$, the geometric complexity measure we use here is based on gaussian averages associated with localizations of the class. We refer the reader to Chapter~12 in \cite{MR1932358} for more details on gaussian processes (in particular to Theorem~12.1.3 for the existence of such a process and to Theorem~12.1.4 for its linearity). 

Denote by $\{G_f : f \in \cF\}$ the canonical gaussian process indexed by $\cF$, that is $\E G_f =0$ and the covariance is given by the inner product in $L_2(\mu)$: $\E G_f G_h = \inr{f, h}_{L_2(\mu)}=\E f(X)h(X)$. Given a set $\cF^\prime \subset \cF$ we put
$$
\E\|G\|_{\cF^\prime}=\sup\Big\{ \E\sup_{h \in \cH} G_h  \ : \ \cH \subset \cF^\prime \ {\rm is \ finite} \Big\}.
$$
This supremum is called the \textit{lattice supremum} (see Chapter~2.2 in \cite{LT:91} for more details). 

As an example, if $\cF^\prime=\{\inr{\cdot, t}:t\in T\}$ is a set of linear functionals indexed  by $T\subset\R^d$ and $X$ is a random vector in $\R^d$ with covariance matrix $\Sigma$ then for $G\sim \cN(0, \Sigma)$, we simply have
\begin{equation*}
\E\|G\|_{\cF^\prime} = \E \sup_{t\in T}\inr{G, t}.
\end{equation*}

We are now in a position to introduce the two complexity parameters that will serve as residual terms in the exact oracle inequalities satisfied by ERM.

\begin{Definition}\label{def:two_fixed_points}
For any $s\geq0$, set $sD=\{f\in L_2(\mu):\norm{f}_{L_2(\mu)}\leq s\}$ and $\cF-\cF=\{f-h:f,h\in\cF\}$. For every $\eta>0$, let
\begin{equation} \label{eq:s-N}
s^*_N(\eta)= \inf \left\{ s>0 : \E\|G\|_{s D \cap ({\cal F}-{\cal F})} \leq \eta s^2 \sqrt{N} \right\},
\end{equation}
and for every $Q>0$, set
\begin{equation} \label{eq:r-N}
r_N^*(Q)= \inf \left\{ r>0 : \E\|G\|_{r D \cap ({\cal F}-{\cal F})} \leq Q r \sqrt{N} \right\}.
\end{equation}
\end{Definition}
In what follows we will always assume without mentioning it explicitly that the sets in \eqref{eq:s-N} and \eqref{eq:r-N} are nonempty (for example, this forces that $Q \geq c/\sqrt{N}$).

There are many situations in which sharp estimates on $\E\|G\|_{r D \cap ({\cal F}-{\cal F})}$ are known and one can identify the fixed points $s^*_N(\eta)$ and $r_N^*(Q)$. We will present several examples of that kind in Section~\ref{sec:examples-TheoB}.

\vskip0.3cm

With these definitions in place, one may formulate a restricted version of the upper bound on the performance of ERM -- for a convex, $L$-subgaussian class of functions.

\vspace{0.5cm}
\noindent {\bf Theorem A.}{\it \hspace{0.2cm}
For every $L \geq 1$ there exist constants $c_1,c_2,c_3$ and $c_4$ that depend only on  $L$ for which the following holds. Let $\cF \subset L_2(\mu)$ be a convex, $L$-subgaussian class of functions, assume that $\|Y-f^*(X)\|_{\psi_2}\leq \sigma$ and set $\eta= c_1/(L\sigma)$ and $Q = c_2/L^2$.
\begin{description}
\item{1.} If $\sigma \geq c_3r_N^*(Q)$ then with probability at least $1-6\exp(-c_4N \eta^2 (s_N^*(\eta))^2)$,
$$
R(\hat{f}) \leq \inf_{f \in \cF} R(f) + (s_N^*(\eta))^2.
$$
\item{2.} If $\sigma \leq c_3r_N^*(Q)$ then with probability at least $1-6\exp(-c_4NQ^2)$,
$$
R(\hat{f}) \leq \inf_{f \in \cF} R(f) + (r_N^*(Q))^2.
$$
\end{description}
Hence, with probability at least
$1-6\exp\left(-c_4N \min\{\eta^2(s_N^*(\eta))^2, Q^2\}\right)$,
$$
R(\hat{f}) \leq \inf_{f \in \cF} R(f) + \max\big\{(s_N^*(\eta))^2,(r_N^*(Q))^2\big\}.
$$}

We will show in what follows that the parameters involved in the upper bound have very clear roles. $r_N^*$ is an upper estimate on the error rate one could have if the problem were noise-free -- that is, if $\sigma=0$. This intrinsic error occurs because it is impossible to distinguish between $f_1,f_2 \in \cF$ using the sample $\bX=(X_i)_{i=1}^N$ when $(f_1(X_i))_{i=1}^N=(f_2(X_i))_{i=1}^N$.

Once noise is introduced to the problem and passes a certain threshold, it is no longer realistic to expect that an intrinsic parameter, which does not depend on the noise level, can serve as an upper bound. And, indeed, $s_N^*(\eta)$ measures the interaction between the `noise'\footnote{We keep the terminology from Statistics: the difference between the output variable $Y$ and the target function $f^*(X)$ is called the noise. This coincides with the classical definition of noise in Statistics when $f^*$ is the regression function.} $f^*(X)-Y$ and the class through the choice of $\eta \sim 1/\sigma$. Thus, beyond a certain noise-level $\sigma$, which depends on the `complexity' of the class $\cF$, $s_N^*(c/\sigma)$ becomes the dominant term in the upper bound.

Note that in the free-noise case, $\sigma=0$, one has $s_N^*(c/\sigma)=0$. Therefore, the error rate of ERM depends only on $r_N^*(Q)$. Also, when the number of observations $N$ is large enough, one also has $r_N^*(Q)=0$, leading to exact reconstruction.

Of course, Theorem~A would be better justified if one could obtain matching lower bounds, showing that ERM is an optimal procedure for subgaussian problems. To that end, it seems natural to employ minimax theory (see, e.g., \cite{MR2724359, yang_barron99,MR1462963,MR1240719,MR804584} for more details on minimax bounds).

What is a reasonable way of identifying a lower bound on the performance of a learning procedure is to see what accuracy and confidence it can guarantee for a minimal set of admissible targets ${\cal Y}$, and a natural choice of a minimal set of targets is 
\begin{equation}\label{eq:gauss-lin-model}
{\cal Y} = \{ Y^f: \  Y^f=f(X)+W\}
\end{equation}
for every $f \in \cF$ and $W$ that is a centered gaussian random variable that has variance $\sigma^2$ and is independent of $X$. Thus, this minimal set of targets consists of `independent perturbations' of realizable learning problems, and thus is arguably the smallest set of `noisy' targets. The minimax rate is (at least) the best accuracy/confidence trade-off that a learning procedure may attain in $\cF$ for the set targets \eqref{eq:gauss-lin-model}. Our main focus will be on the accuracy/confidence tradeoff for the accuracy level described in Theorem~A.

Standard minimax bounds are based on information-theoretical results such as Fano's Lemma, Assouad's Lemma or Pinsker's inequalities. Unfortunately, these results do not yield lower bounds in the high probability realm of Theorem A; rather, these results are restricted to constant confidence or hold in expectation. To treat the high probability regime, we present a new minimax bound that is based on the gaussian shift theorem (and therefore on the gaussian isoperimetric inequality).

\vspace{0.6cm}
\noindent{\bf Theorem A${}^\prime$.}{\it \hspace{0.2cm} There exists an absolute constant $c_5$ for which the following holds. Let $\cF \subset L_2(\mu)$ be a class that is star-shaped around one of its points (i.e., for some $f_0\in\cF$ and every $f\in\cF$, $[f_0,f] \subset \cF$), and let ${\cal Y}$ be the set of admissible targets from \eqref{eq:gauss-lin-model}. If $\tilde f_N$ attains an accuracy $\eps_N$ with the confidence level $\delta_N$ for any target $Y^f \in {\cal Y}$, then
$$
\eps_N \geq \min \left\{c_5\sigma^2 \frac{\log(1/\delta_N)}{N}, \frac{1}{4}{\rm diam}^2(\cF,L_2(\mu))\right\}.
$$
}
\vskip0.3cm

Note that no assumption on the underlying measure $\mu$ is required in Theorem A${}^\prime$. Moreover, Theorem A${}^\prime$ makes a natural connection between accuracy and confidence: the higher the confidence $1-\delta_N$ the larger $\eps_N$ must be. 

An important outcome of Theorem~A and Theorem~A${}^\prime$ is that for the set of admissible targets ${\cal Y}$ as in \eqref{eq:gauss-lin-model}, and as long as the class $\cF$ is convex and $L$-subgaussian, ERM is optimal in the following sense:

\vspace{0.6cm}
\noindent{\bf Theorem A${}^{\prime\prime}$.}{\it \hspace{0.2cm}
There exist absolute constants $c_1,...,c_4$ for which the following holds. Let $\cF$ be a convex, $L$-subgaussian class of functions and consider the set of admissible targets ${\cal Y}$ as in \eqref{eq:gauss-lin-model}. Set $\eta = c_1/(L\sigma)$ and $Q =  c_2/L^2$. If $\sigma \geq c_3r_N^*(Q)$ then for any target $Y^f \in {\cal Y}$, the ERM $\hat f$ satisfies
$$
R(\hat{f}) \leq \inf_{f \in \cF} R(f) + (s_N^*(\eta))^2 \ \  {\rm with \ probability} \ \ 1-6\exp\big(-c_4N\eta^2 (s_N^*(\eta))^2\big).
$$
Also, for any learning procedure $\tilde{f}$ there is some $f \in \cF$ for which, if given the data generated by the target $Y^f$ and $R(\tilde{f}) \leq \inf_{f \in \cF} R(f) + (s_N^*(\eta))^2$ with probability at least $1-\delta$, then
$$
\delta \geq  \exp\big(-c_5N\eta^2 (s_N^*(\eta))^2\big).
$$
}
Thus, up to the constant in the exponent, the upper bound and the lower bound match and the ERM achieves this bound.
\vskip0.5cm

The second question we wish to address is what happens when the desired confidence is an absolute constant -- for example, when $1-\delta_N$ is, say, $1/2$, but the noise level is nontrivial in the sense that $s_N^*$ dominates $r_N^*$. We will show that in such a situation, Theorem~A is optimal in a minimax sense under some regularity assumptions on $\cF$. This complements Theorem A${}^{\prime\prime}$ which proves the optimality of ERM (under no extra structural assumption) in the high probability case -- when $\delta_N\sim\exp(-c\eta^2 (s_N^{*}(\eta))^2 N)$.

To explore the constant confidence regime, let us consider the `Sudakov analog' of the gaussian-based parameter $s_N^*(\eta)$: recall that by Sudakov's inequality (see, for example, \cite{LT:91}), for any $r>0$,
\begin{equation}\label{eq:sudakov-ineq}
\sup_{\eps >0 } \eps \log^{1/2} N((\cF -\cF) \cap rD, \eps D) \lesssim \E\|G\|_{rD \cap (\cF -\cF)}.
\end{equation}
Put $C(r)= \sup_{f\in\cF}r \log^{1/2} N((\cF -f) \cap 2rD, r D)$ and set
$$
q_N^*(\eta) = \inf \{s>0 : C(s) \leq \eta s^2 \sqrt{N} \}.
$$

\vskip0.5cm
\noindent {\bf Theorem B.}{\it \hspace{0.2cm}
There exists an absolute constant $c_1$ for which the following holds. Let $\cF$ be a class of functions, set $W \sim \cN(0,\sigma^2)$ and for every $f \in \cF$, put $Y^f=f(X)+W$. If $\tilde f_N$ performs with a confidence parameter $\delta_N<1/4$ for every such target $Y^f$,  then its accuracy cannot be better than $ c_1 (q_N^*(c_1/\sigma))^2$.
}
\vskip0.5cm

Theorem B is known, and may be derived from Theorem~2.5 in \cite{MR2724359} or from \cite{yang_barron99}. The proof presented here is new, and follows the same path as the proof of Theorem A${}^{\prime}$.

With Theorem A in mind, Theorem B implies that if the learning problem is subgaussian, $s_N^*(\eta)$ and  $q_N^*(\eta^\prime)$ are equivalent for $\eta, \eta^\prime \sim 1/\sigma$ and $\sigma\gtrsim r_N^*$, the minimax rate in the constant probability regime is attained by ERM.

Finally, let us consider the low-noise case, in which $\sigma \lesssim r_N^*$. Although it is not clear if $r_N^*$ is an optimal bound in that range (except when $\sigma \sim r_N^*$), it turns out that it is not far from optimal.

\begin{Definition} \label{def:gen-Gelfand}
Let $\cF$ be a class of functions. For every sample $\bX=(X_1,...,X_N)$ and $f \in \cF$, set
$$
K(f,\bX) = \{ h \in \cF : (f(X_i))_{i=1}^N = (h(X_i))_{i=1}^N\},
$$
which is the ``level set" in $\cF$ given by the values of $f$ on the sample. Let $\cD(f,\bX)$ be the diameter of $K(f,\bX)$ with respect to the $L_2(\mu)$ norm.
\end{Definition}

Clearly, if $\sigma =0$ then for every sample $\bX$, ERM selects $\hat{f} \in K(f^*,\bX)$ and since $Y=f^*(X)$, $R(f)=\|f-f^*\|_{L_2(\mu)}^2$. Thus, $R(\hat{f}) \leq \cD^2(f^*,\bX)$. It is natural to ask whether the reverse direction is true. The following result shows that the largest typical value of $\cD(f,\bX)$ is a constant-probability  minimax bound.
\vskip0.6cm
\noindent {\bf Theorem~C.}{\it \hspace{0.2cm}
For every $f \in \cF$ and $V$ that is independent of $X$, set $Y^f=f(X)+V$. Then, for any procedure $\tilde{f}_N$,
\begin{equation*}
  \sup_{f \in \cF} Pr\left(\|\tilde f_N((Y^f_i,X_i)_{i=1}^N)-f\|_{L_2(\mu)} \geq  \frac{1}{4}\cD(f,\bX) \right) \geq 1/2,
\end{equation*}
with the probability taken with respect to the product measures endowed by $(Y_i^f,X_i)_{i=1}^N$.
}
\vskip0.5cm

One natural example in which Theorem~C may be used is when $T$ is a convex, centrally-symmetric subset of $\R^d$ (i.e., if $t \in T$ then $-t \in T$), and $\cF$ is the class of linear functionals indexed by $T$, i.e., $\cF=\{\inr{t,\cdot} : t \in T\}$. Let $X_1,...,X_N$ be an independent sample selected according to an isotropic probability measure on $\R^d$. If $\{e_1,\ldots,e_N\}$ is the canonical basis of $\R^N$ and $\Gamma=\sum_{i=1}^N \inr{X_i,\cdot}e_i$ is the random matrix whose rows are $(X_i)_{i=1}^N$, then $\cD(0,\bX)$ is the diameter of the intersection of the kernel of $\Gamma$ and $T$ : $\cD(0,\bX)= {\rm ker}(\Gamma) \cap T$. If $f^*=\inr{t^*,\cdot}$, one can relate  $\cD(t^*,\bX)$ to the Gelfand widths of $T$ (see, e.g., \cite{MR1036275} and \cite{MR774404} for more details).

\begin{Definition}\label{def:gelfand_width}
Let $T$ be a convex, centrally-symmetric subset of $\R^d$. The \textbf{Gelfand $N$-width of $T$} is the smallest $\ell_2^d$-diameter of an $N$-codimensional section of $T$. In particular,
\begin{equation*}
c_N(T) = \inf\left\{ {\rm diam}({\rm ker}(\Gamma) \cap T, \ell_2^d) : \Gamma\in L(\R^d, \R^N)\right\},
\end{equation*}
where $L(\R^d, \R^N)$ is the set of all linear operators from $\R^d$ to $\R^N$ and ${\rm diam}(V, \ell_2^d) = \sup_{u,v\in V}\norm{u-v}_2$.
\end{Definition}

Hence, for every $t_0 \in T$,
$$
c_N(T) \leq {\rm diam}\left(K(t_0,\bX)-t_0,\ell_2^d\right) \leq 2\cD(0,\bX),
$$
and by Theorem C,  $c_N(T)/8$ is a lower bound on the minimax rate in the constant confidence regime. Therefore, when $r_N^* \sim c_N(T)$, it follows that for every $0 \leq \sigma \lesssim r_N^*$, $r_N^*$ is the constant-probability minimax rate, and that rate is achieved by ERM.

\vskip0.3cm

It should be noted that although our presentation focuses on oracle inequalities in a given class, oracle inequalities for model selection and regularized procedures can be derived from the isomorphic method in general, specifically, from Theorem~\ref{thm:ratio-estimate} below. This strategy is rather standard and has been used, for example, in \cite{MR2981422,MR2590050}, in Chapter~3.6  of \cite{HDR} or recently in \cite{LM_sparsity}. We will not present results on regularization methods or model selection methods in what follows since those may be easily obtained from results on ERM.

\vskip0.3cm

We end this introduction with a word about notation. Throughout, absolute
constants or constants that depend on other parameters are denoted by $c$, $C$, $c_1$,
$c_2$, etc., (and, of course, we will specify when a constant is
absolute and when it depends on other parameters); their values may change from line to line. The notation $x\sim y$ (resp. $x\lesssim y$) means that there exist absolute constants $0<c<C$ for which $cy\leq x\leq Cy$ (resp. $x\leq Cy$). If $b>0$ is a parameter then $x\lesssim_b y$ means that $x\leq C(b) y $ for some constant $C(b)$ that depends only on $b$.

Let $\ell_p^d$ be $\R^d$ endowed with the norm $\norm{x}_{\ell_p^d}=\big(\sum_{j=1}^d |x_j|^p\big)^{1/p}$. The unit ball in $\ell_p^d$ is denoted by $B_p^d$, and the unit Euclidean sphere in $\R^d$ is $S^{d-1}$. We also denote by $d_{L_2}(\cF^\prime)$ the  diameter of $\cF^\prime$ in $L_{2}(\mu)$.

The proofs of our main results are presented in the next two sections. We then present several examples of applications of those results, in which the rates established in Theorem A are shown to be sharp in both the high and constant confidence regimes. The final section contains some concluding remarks.

\section{Proof of Theorem A} \label{sec:proof-TheoB+examples}
The proof of Theorem~A shows that it is more general than stated. Rather than convexity, the two properties that are actually needed are the following:
\begin{Definition} \label{def:star-shape}
A class $\cH$ is star-shaped around $h_0 \in \cH$ if for every $h \in \cH$, the interval $[h,h_0]$ is contained in $\cH$.
\end{Definition}
We will assume that $\cF-\cF = \{f-h: f,h \in \cF\}$ is star-shaped around $0$, otherwise, one may consider the star-shaped hull of $\cF-\cF$ with $0$, that is, the set
$$
\left\{\lambda(f-h) : 0 \leq \lambda \leq 1, \ f,h \in \cF \right\}
$$
which is not much larger than $\cF-\cF$.
\vskip0.3cm
The second property required is a variant of the Bernstein condition (cf. \cite{MR2240689}).
\begin{Definition} \label{def:Bernstein-star-shape}
A class $\cF$ is \textit{$B$-Bernstein} relative to the target $Y$, if for every $f\in\cF$,
  \begin{equation}\label{def:one-sided-Bernstein}
    \E\big(f(X)-f^*(X)\big)^2\leq BP\cL_f = B \E\big((Y-f(X))^2-(Y-f^*(X))^2\big).
  \end{equation}
\end{Definition}
Definition \ref{def:Bernstein-star-shape} is far less restrictive than it appears at first glance. Indeed, by the $2$-convexity of the $L_2$ norm, if $\cF$ is convex then for any target $Y \in L_2$, $\cF$ is $1$-Bernstein relative to $Y$. Moreover, the results from \cite{MR2426759} show that for every class $\cF$ and every target $Y$, the Bernstein constant depends only on the distance between $Y$ and the set of targets $Z$ for which the functional $f\to\E(f-Z)^2$ has multiple minimizers in $\cF$. Finally, note that if one wishes $\cF$ to satisfy a Bernstein condition relative to {\it every} target $Y$, it forces $\cF$ to be convex in the locally-compact case (see Section \ref{sec:conc-rem} for more details).

In what follows, we shall assume that $\cF-\cF$ is star-shaped around $0$ and that $\cF$ satisfies the Bernstein condition \eqref{def:one-sided-Bernstein}.
\vskip0.3cm

The next lemma (which will be proved in the Appendix) shows that the assumption that $\cF-\cF$ is star-shaped around $0$ adds some regularity to the gaussian process $\{G_f : f \in \cF-\cF\}$.

\begin{Lemma}\label{lem:complexity-regularity}
Assume that $\cF-\cF$ is star-shaped around $0$ and let $\psi:s\geq0\to\E\norm{G}_{sD\cap(\cF-\cF)}$. Then the following holds:
\begin{enumerate}
  \item $\phi:s\to\psi(s)/s$ is non-increasing.
  \item For $\eta>0$ and any $s\geq s_N^*(\eta)$, $\psi(s)\leq \eta s^2 \sqrt{N}$,  and for any $0<s<s_N^*(\eta)$, $\psi(s) \geq \eta s^2 \sqrt{N}$.
   \item Let $Q>\sqrt{\pi/2N}$. For any $r\geq r_N^*(Q)$, $\psi(r)\leq Q r \sqrt{N}$ and for any $0<r<r_N^*(Q)$, $\psi(r) > Q r\sqrt{N}$.
\end{enumerate}
\end{Lemma}

A straightforward outcome of Lemma \ref{lem:complexity-regularity} which will be used later is as follows:
\begin{Lemma} \label{lemma:r-vs-s}
Let $c,\sigma,Q >0$, set $\eta=c/\sigma$ and consider $s_N^*(\eta)$ and $r_N^*(Q)$ as introduced in Definition~\ref{def:two_fixed_points}.
\begin{enumerate}
\item If $\sigma \geq (c/Q) r_N^*(Q)$ then $s_N^*(\eta) \geq r_N^*(Q)$, and if $\sigma \leq (c/Q)r_N^*(Q)$ then $s_N^*(\eta) \leq r_N^*(Q)$.
\item If $s_N^*(\eta) \geq r_N^*(Q)$ then $ \eta s_N^*(\eta) \leq 4 Q$.
\end{enumerate}
\end{Lemma}
The proof of Lemma \ref{lemma:r-vs-s} will also be presented in the Appendix.

\vskip0.3cm

When considering the parameters $r_N^*(Q)$ and $s_N^*(\eta)$, what may seem odd at first glance is the different normalization in their definition -- the first condition is linear, while the second is quadratic. The two originate from the need to compare the way in which two processes, the quadratic component and the multiplier component of the excess loss functional scale with $\|f-f^*\|_{L_2(\mu)}$. Indeed, note that
\begin{align*}
{\cal L}_f(X,Y) =  (f(X)-Y)^2-(f^*(X)-Y)^2
=  (f(X)-f^*(X))^2+2(f(X)-f^*(X))(f^*(X)-Y).
\end{align*}
The quadratic term (i.e. $(f(X)-f^*(X))^2$) is noise-free, and as will be explained below, $r_N^*$ measures the lowest level $r$ at which if $\|f-f^*\|_{L_2(\mu)} \geq r$, then $\E(f-f^*)^2 \sim N^{-1}\sum_{i=1}^N (f-f^*)^2(X_i)$.

In contrast, $s_N^*$ is designed for dealing with the multiplier process, originating from the term $(f^*(X)-Y)\cdot(f-f^*)(X)$. To compare the resulting multiplier component with $\E(f-f^*)^2$ (which is the order of magnitude of $N^{-1}\sum_{i=1}^N (f-f^*)^2(X_i)$ when $\|f-f^*\|_{L_2(\mu)} \geq r_N^*$), one has to study
$$
f \to \frac{1}{N}\sum_{i=1}^N (f^*(X_i)-Y_i) \cdot \frac{(f-f^*)(X_i)}{\E (f-f^*)^2},
$$
and that is the source of the seemingly less-natural normalization in the definition of $s_N^*(\eta)$.

\vskip0.3cm

Let us begin with an estimate on the quadratic component, which is based on a functional Bernstein type inequality (see \cite{MR3354613, bednorz, shahar_proc}).

\begin{Theorem}\label{thm:quadratic}
There exist absolute constants $c_1$ and $c_2$ for which the following holds.
Let ${\cal H}$ be an $L$-subgaussian class. For every $u>0$, with probability at least $1-2\exp(-c_1\min(u^2, u\sqrt N))$,
\begin{equation}
    \label{eq:proc_quad}
    \sup_{h\in \cH}
    \left|\frac{1}{N}\sum_{i=1}^N h^2(X_i)-\E h^2\right|
    \leq c_2 L^2\left(
        \frac{ d \gamma}{\sqrt{N}}
        + \frac{\gamma^2}{N}
        + \frac{u d^2}{\sqrt N}
    \right)
\end{equation}
where $d=d_{L_2}(\cH)$ is the diameter in $L_2(\mu)$ of $\cH$ and $\gamma = \E\norm{G}_{\cH}$.
\end{Theorem}

\vskip0.3cm

The following result is a straightforward application of Theorem \ref{thm:quadratic} and illustrates the role of $r_N^*(Q)$.
\begin{Lemma} \label{lemma:above-r-N}
There exist absolute constants $c_1,c_2$ and $c_3$ for which the following holds. Let $\cF$ be an $L$-subgaussian class, assume that $\cF-\cF$ is star-shaped around $0$ and let $f^* \in \cF$. If $0<Q \leq 1$ and $r > r_N^*(Q)$, then with probability at least
$1-2\exp\big(-c_1Q^2 N\big)$,
\begin{equation*}
\sup_{h \in rD \cap (\cF -f^*)}  \left|\frac{1}{N}\sum_{i=1}^N h^2(X_i)-\E h^2 \right| \leq c_2QL^2 r^2.
\end{equation*}
\end{Lemma}
\proof The claim is an immediate corollary of Theorem~\ref{thm:quadratic}. Indeed, one simply has to apply Theorem~\ref{thm:quadratic} to the set $\cH=rD\cap(\cF-\cF)$ and to recall that by Lemma~\ref{lem:complexity-regularity}, if $r > r_N^*(Q)$ then $\E\|G\|_{rD \cap (\cF-\cF)}\leq Qr\sqrt{N}$. Therefore, for any $u>0$, with probability at least $1-2\exp(-c_1\min(u^2, u\sqrt N))$,
\begin{equation*}
    \sup_{f,h\in \cF:\norm{f-h}_{L_2(\mu)}\leq r}
    \left|\frac{1}{N}\sum_{i=1}^N (f-h)^2(X_i)-\E (f-h)^2\right|
    \leq c_2 L^2\left(
        \frac{ d \gamma}{\sqrt{N}}
        + \frac{\gamma^2}{N}
        + \frac{u d^2}{\sqrt N}
    \right)
\end{equation*}
where $d={\rm diam}(rD\cap(\cF-\cF), L_2)\leq r$ and $\gamma = \E\|G\|_{rD \cap (\cF-\cF)}\leq Qr\sqrt{N}$. Hence, for $u=c_2 Q \sqrt{N}$, with probability larger than $1-2\exp(-c_3 Q^2 N)$,
\begin{equation}\label{eq:RIP-shahar-2}
    \sup_{f,h\in \cF:\norm{f-h}_{L_2(\mu)}\leq r}
    \left|\frac{1}{N}\sum_{i=1}^N (f-h)^2(X_i)-\E (f-h)^2\right|
    \leq c_4 Q L^2 r^2.
\end{equation}
\endproof

\begin{remark}
Using the notation of Lemma~\ref{thm:quadratic}, consider $Q \leq \min\{1/(2c_2L^2),1\}$. If \eqref{eq:RIP-shahar-2} holds then for every $f \in \cF$ that satisfies $\|f-f^*\|_{L_2(\mu)} \geq r > r_N^*(Q)$, one clearly has
$$
\frac{1}{2} \E(f-f^*)^2 \leq \frac{1}{N}\sum_{i=1}^N (f-f^*)^2(X_i) \leq \frac{3}{2}\E(f-f^*)^2;
$$
this is evident because for $h=f-f^*$,
$$
\left|\frac{1}{N}\sum_{i=1}^N h^2(X_i) - \E h^2 \right| \leq \frac{r^2}{2} \leq \frac{\E h^2}{2}.
$$
\end{remark}

The second ingredient required for the proof of Theorem A is a bound on multiplier processes.
\begin{Theorem} \label{thm:multiplier}[Theorem~4.4 in \cite{shahar_proc}]
There exist absolute constants $c_1$ and $c_2$ for which the following holds. If ${\cal H}$ is an $L$-subgaussian class and $\xi \in L_{\psi_2}$, then for every $u,w \geq 8$, and every integer $s_0\geq1$,
with probability at least
\begin{equation*}
1-2\exp\big(-c_1 u^2 2^{s_0}) - 2 \exp\big(-c_1 N w^2\big),
\end{equation*}
\begin{equation*}
\sup_{h \in {\cal H}} \left|\frac{1}{N}\sum_{i=1}^N \xi_i h(X_i) - \E\xi h(X) \right| \leq c_2L u w \frac{\|\xi\|_{L_{\psi_2}}}{\sqrt{N}}\left(\E\norm{G}_\cH + 2^{s_0/2}d_{L_2}(\cH)\right).
\end{equation*}
\end{Theorem}

Note that in Theorem \ref{thm:multiplier} one does not assume that $\xi$ and $X$ are independent, a fact that will be significant in what follows. Indeed, we will apply Theorem~\ref{thm:multiplier} to  $\xi = Y-f^*(X)$ and the class $\cH = r D \cap(\cF-\cF)$ for $r > s_N^*(\eta)$. In that case, $d_{L_2}(\cH)\leq r$ and $\E\norm{G}_\cH\leq \eta r^2 \sqrt{N}$, and for $2^{s_0/2}\sim \eta r \sqrt{N}$, we obtain that with probability larger than  $1-4\exp\left(-c_1 N \min\{\eta^2 r^2, 1\}\right)$,
\begin{equation}\label{eq:coro_shahar_multipl}
\sup_{f,h\in\cF:\norm{f-h}_{L_2(\mu)}\leq r} \left|\frac{1}{N}\sum_{i=1}^N \xi_i (f-h)(X_i) - \E\xi (f-h)(X) \right| \leq c_2L \eta \|\xi\|_{L_{\psi_2}} r^2.
\end{equation}

Combining the estimates on the quadratic and multiplier process leads to the following ratio estimate:
\begin{Theorem}\label{thm:ratio-estimate}
For every $L \geq 1$ and  $B\geq1$ there exist constants $c_0,c_1,c_2$ and $c_3$ that depend only on $B$ and $L$ for which the following holds. Let ${\cal F}$ be an $L$-subgaussian class that is $B$-Bernstein relative to the target $Y$. Assume that $\cF-\cF$ is star-shaped around $0$ and that $\|Y-f^*(X)\|_{\psi_2} \leq \sigma$. Set $\eta = c_0/(LB\sigma)$ and $Q=c_1/(L^2 B)$.
\begin{description}
\item{1.} If $\sigma \geq c_2r_N^*(Q)$, then with probability at least
$1-6\exp\left(-c_3N \cdot \eta^2(s_N^*(\eta))^2\right)$,
$$
\sup_{\{f\in\cF : P\cL_f \geq (s_N^*(\eta))^2/B\}} \left|\frac{1}{N} \sum_{i=1}^N \frac{{\cal L}_f(X_i,Y_i)}{P {\cal L}_f} -1 \right| \leq \frac{1}{2}.
$$
\item{2.} If $\sigma \leq c_2r_N^*(Q)$, then with probability at least
$1-6\exp\left(-c_3Q^2N/B\right)$,
$$
 \sup_{\{f\in\cF : P\cL_f \geq (r_N^*(Q))^2/B\}} \left|\frac{1}{N} \sum_{i=1}^N \frac{{\cal L}_f(X_i,Y_i)}{P {\cal L}_f} -1 \right| \leq \frac{1}{2}.
$$
\end{description}
\end{Theorem}

\proof
Set $\xi=(f^*(X)-Y)$ and thus
\begin{align*}
{\cal L}_f(X,Y) & = (f-f^*)^2(X)+2\xi(f-f^*)(X).
\end{align*}
Fix $\lambda>0$ and let ${\cal F}_\lambda =\{ f\in\cF : P {\cal L}_f \geq \lambda \}$. Since $\cF$ satisfies the $B$-Bernstein condition relative to $Y$, it follows that for every $f\in\cF$, $\|f-f^*\|^2_{L_2(\mu)} \leq B P{\cal L}_f$. Moreover, if $f \in {\cal F}_\lambda$ then
\begin{equation} \label{eq:subgaussian-ratio-1}
\left\|\frac{f-f^*}{(P{\cal L}_f)^{1/2}}\right\|_{L_2(\mu)}^2 \leq B \ \
\mbox{ and } \ \
\left\|\frac{f-f^*}{P{\cal L}_f}\right\|_{L_2(\mu)}^2 \leq \frac{B}{P{\cal L}_f} \leq \frac{B}{\lambda}.
\end{equation}
Therefore,
\begin{align*}
\sup_{f \in {\cal F}_\lambda} & \left|\frac{1}{N}\sum_{i=1}^N \frac{{\cal L}_f(X_i,Y_i)}{P {\cal L}_f}-1 \right| =   \sup_{f \in {\cal F}_\lambda} \left|\frac{1}{N}\sum_{i=1}^N \frac{{\cal L}_f(X_i,Y_i)-P {\cal L}_f}{P {\cal L}_f} \right|
\\
\leq & \sup_{f \in {\cal F}_\lambda}  \left|\frac{1}{N}\sum_{i=1}^N \left(\frac{f-f^*}{\big(P{\cal L}_f\big)^{1/2}}\right)^2(X_i)-\E\left(\frac{f-f^*}{(P{\cal L}_f)^{1/2}}\right)^2 \right|
+
2\sup_{f \in {\cal F}_\lambda } \left|\frac{1}{N}\sum_{i=1}^N \xi_i \left(\frac{f-f^*}{P{\cal L}_f}\right)(X_i) - \frac{\E \xi (f-f^*)}{P\cL_f} \right|.
\end{align*}
Set
$$
W_\lambda = \left\{\frac{f-f^*}{(P{\cal L}_f)^{1/2}} : f \in {\cal F}_\lambda \right\},  \  \ V_\lambda = \left\{\frac{f-f^*}{P{\cal L}_f} : f \in {\cal F}_\lambda \right\},
$$
and ${\cal H}=({\cal F}-{\cal F}) \cap \sqrt{\lambda B}D$. Recall that $\cF-\cF$ is star-shaped around $0$, and by (\ref{eq:subgaussian-ratio-1}) one has that
$$
W_\lambda \subset \frac{1}{\sqrt{\lambda}}({\cal F}-{\cal F}) \cap \sqrt{B} D \subset \frac{1}{\sqrt{\lambda}}\left(({\cal F}-{\cal F}) \cap \sqrt{\lambda B} D\right)=\frac{{\cal H}}{\sqrt{\lambda}},
$$
and
$$
V_\lambda \subset \frac{1}{\lambda}({\cal F}-{\cal F}) \cap \left(\sqrt{\frac{B}{\lambda}}\right)D \subset \frac{1}{\lambda} \left(({\cal F}-{\cal F}) \cap \sqrt{\lambda B}D \right)=\frac{{\cal H}}{\lambda}.
$$

Fix  $\eta = c_0/(LB\sigma)$ and $Q=c_1/(L^2B)$ for suitable absolute constants $c_0$ and $c_1$. Set $r > r_N^*(Q)$ and note that by Lemma \ref{lemma:r-vs-s}, if $\sigma \geq c_2r_N^*(Q)$  then $r_N^*(Q) \leq s_N^*(\eta)$ and $\eta s_N^*(\eta) \geq 4Q$, and if $\sigma \leq c_2r_N^*(Q)$ then $r_N^*(Q) \geq s_N^*(\eta)$; also $c_2=c_0/LBQ = c_0L/c_1$.

First, consider the case $\sigma \geq c_2r_N^*(Q)$. Applying Lemma~\ref{lemma:above-r-N} for $\lambda=(s_N^*(\eta))^2/B$, it follows that with probability at least $1-2\exp(-c_3 Q^2N)$
$$
\sup_{w \in W_\lambda} \left|\frac{1}{N}\sum_{i=1}^N w^2(X_i)-\E w^2\right| \leq c_4 Q L^2 B \leq \frac{1}{4},
$$
provided that $Q\leq 1/(4 c_4 L^2 B)$. Moreover, by \eqref{eq:coro_shahar_multipl}, and because $\eta s_N^*(\eta) \geq 4Q$, one has that with probability at least
$$
1-4\exp(-c_5 N \eta^2 (s_N^*(\eta))^2),
$$
\begin{equation*}
\sup_{v \in V_\lambda} \left|\frac{1}{N}\sum_{i=1}^N \xi_i v(X_i) -\E \xi v\right| \leq c_6 L B \sigma \eta  \leq \frac{1}{8}
\end{equation*}
as long as $\eta\leq 1/(c_7 L B \sigma)$.

Thus, for any $Q \lesssim 1/L^2B$ and $\eta \lesssim 1/\sigma L B$, if $\sigma \geq c_2 L r_N^*(Q)$ then with probability at least $1-6\exp(-c_8N \cdot \eta^2(s_N^*(\eta))^2)$, the following holds: for every $f\in\cF$ that satisfies that $P\cL_f \geq \lambda$,
$$
\frac{1}{2} P\cL_f \leq P_N \cL_f \leq \frac{3}{2}P\cL_f.
$$

Next, let us consider that case $\sigma \leq c_2 r_N^*(Q)$ which follows a very similar path to the first case. Recall that $r_N^*(Q)\geq s_N^*(\eta)$. Setting $\lambda=(r_N^*(Q))^2/B$, it follows from Lemma \ref{lemma:above-r-N} and \eqref{eq:coro_shahar_multipl} that with probability at least
\begin{equation}\label{eq:proba_1}
1-2\exp(-cQ^2N) - 4\exp\left(-c N\min\Big\{\frac{\eta^2(r_N^*(Q))^2}{\sigma^2 B}, 1\Big\}\right),
\end{equation}
\begin{equation*}
\sup_{w \in W_\lambda} \left|\frac{1}{N}\sum_{i=1}^N w^2(X_i) - \E w^2 \right| \leq \frac{1}{4} \mbox{ and }
\sup_{v \in V_\lambda} \left|\frac{1}{N}\sum_{i=1}^N \xi_iv(X_i) - \E \xi v \right| \leq \frac{1}{8}
\end{equation*}
as long as  $Q \lesssim 1/(L^2B)$ and $\eta \lesssim 1/(\sigma L B)$. The claim now follows because $\eta \sim 1/\sigma$ and by the choice of $\sigma$, namely, that $r_N^*(Q)/\sigma \geq c_2$.
\endproof

Theorem A is an immediate outcome of Theorem \ref{thm:ratio-estimate} for $B=1$ and the isomorphic method described in the introduction.

\section{Minimax lower bounds (proofs of Theorem~A${}^{\prime}$, B and C)}
\label{sec:Proof-TheoD}
Let $\cF$ be a class of functions on a probability space $(\Omega,\mu)$, fix $f \in \cF$, let $W$ be a centred gaussian random variable that is independent of $X$ and consider the target function $Y^f=f(X)+W$. For any $\bX=(x_1,\ldots, x_N)\in \Omega^N$, let $\nu_{f,\bX}$ be the conditional probability measure of $(Y^f_i|X_i=x_i)_{i=1}^N$, which is given by
\begin{equation*}
d\nu_{f,\bX}(y)=\exp\left(-\frac{\|y-(f(x_i))_{i=1}^N\|_{\ell_2^N}^2}{2\sigma^2}\right) \cdot \frac{dy}{(\sqrt{2\pi}\sigma)^N},
\end{equation*}
and set $\nu_{f,\bX} \otimes \mu^N $ to be the probability measure on $(\R \times \Omega)^N$ that generates the sample $(Y^f_i,X_i)_{i=1}^N$.

Let
$$
{\cal B}(f,r) = \{h \in \cF : \E{\cal L}_h \leq r\}=\{h \in \cF : \E(f-h)^2 \leq r\},
$$
where $\cL_h(X,Y^f)=(Y^f-h(X))^2-(Y^f-f(X))^2$.

If a procedure $\tilde f_N$ performs with accuracy $\eps_N$ and has a confidence parameter $\delta_N$, then for every $f \in \cF$,
$$
(\nu_{f,\bX}\otimes\mu^N) \left(\tilde f_N^{-1}({\cal B}(f,\eps_N))\right) \geq 1-\delta_N.
$$
In other words, for every $f \in \cF$, the set of data points $(y_i,x_i)_{i=1}^N$ that are mapped by the procedure $\tilde f_N$ into the set $\{h \in \cF: \E {\cal L}_h \leq \eps_N\}$ is of $\nu_{f,\bX}\otimes \mu^N$ measure  at least $1-\delta_N$.

\vskip0.3cm
The first estimate presented here is the high probability lower bound, formulated in Theorem A$^\prime$.

\begin{Theorem} \label{thm:lower-bound-general}
There exists an absolute constant $c_1$ for which the following holds. If $\cF$ is star-shaped around one of its points and $\tilde f_N$ is a procedure that performs with accuracy $\eps_N$ for any target of the form $Y^f$ with a confidence parameter $\delta_N < 1/4$, then
$$
\eps_N \geq \min \left\{c_1 \sigma^2 \frac{\log(1/\delta_N)}{N}, \frac{1}{4}(d_{\cF}(L_2))^2\right\}.
$$
\end{Theorem}

\begin{figure}[h!]
  \begin{tikzpicture}[scale=0.65, line width=3pt]
    \draw [->] (-0.05,-2) -- (12,-2) node[anchor=north] {$\cX^N$}; \draw
    [->] (0,-2.05) -- (0,10) node[anchor=east] {$\R^N$}; \draw (20,5)
    circle (5cm); \draw (20.5,-0.5) node[anchor=east] {$\cF$};
    \filldraw [blue] (7,9) circle (1pt) node[anchor=north]
      {$\cD=(x_i,y_i)_{i=1}^N$}; \filldraw [blue] (16,7) circle (1pt)
      node[anchor=west] {$\tilde{f}_N(\cD,\cdot)$}; \draw[->,thick,
      blue] (7,9) .. controls (10,10) and (15,8) .. (16,7);
    \filldraw [thick] (21,7) circle (2pt) node[anchor=west]
      {$f^*_1$}; \draw [thick] (21,7) circle (2cm);
      \draw
      (23,5.3) node[anchor=north] {$\cB(f^*_1,\eps_N)$};
     \draw[<-,thick, blue] (7,5) .. controls (12,5)
      and (17,5) .. (20,7); \draw[thick] (6,4) ellipse (4cm and 2cm);
        \draw (0.7,5) node[anchor=west] {$\cA_{f^*_1}$};\draw[blue] (13,5) node[anchor=south] {$(\tilde f_N)^{-1}$};

      \draw (5,-2) node[anchor=north] {$\bX=(x_i)_{i=1}^N$};
      \draw[dashed, thick] (5,-2) -- (5,10);
      \draw (5,2) -- (5,6); \draw (5,3) node[anchor=west] {$\cA_{f_1^*}|\bX$};
    \filldraw [thick] (21,7) circle (2pt) node[anchor=west]
      {$f^*_1$};\draw [thick] (21,7) circle (2cm); \filldraw [thick] (20,2.5) circle (2pt) node[anchor=west]
      {$f^*_2$};\draw [thick] (20,2.5) circle (2cm);
\draw[<-,thick, blue] (7,1) .. controls (12,1)
      and (17,1) .. (20,3); \draw[thick] (6,0) ellipse (4cm and 1.5cm);
        \draw (0.7,1) node[anchor=west] {$\cA_{f^*_2}$}; \draw (5,-0.5) node[anchor=west] {$\cA_{f_2^*}|\bX$}; \draw (5,-1.5) -- (5,1.5);
  \end{tikzpicture}
          \caption{Proof of the minimax lower bounds via the gaussian shift theorem in $\R^N$}
        \label{fig:minimax_lower_bound}
\end{figure}
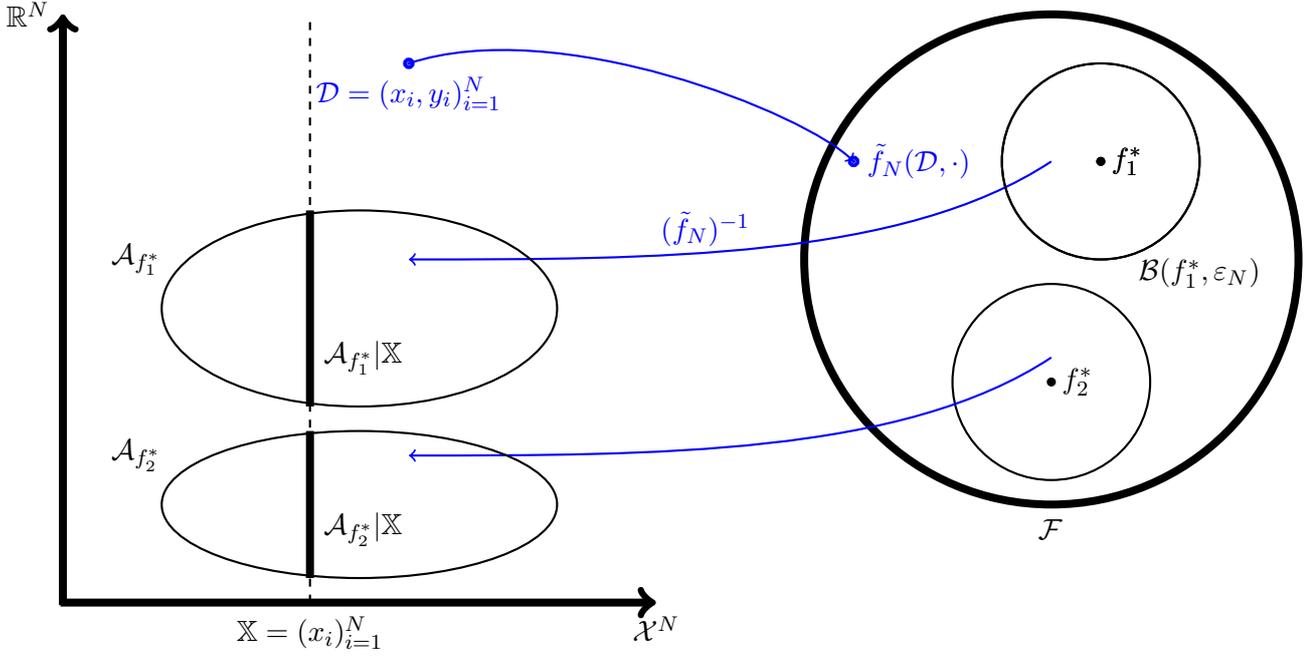

Theorem \ref{thm:lower-bound-general} leads to the lower estimate in Theorem A${}^{\prime\prime}$. Indeed, if a procedure performs with confidence $\delta_N =\exp(-c_0\gamma N)$ for some $\gamma$, then $\eps_N \geq c_2\sigma^2 \gamma$. Setting $\gamma=c_3 \eta(s_N^*(\eta))^2$ for $\eta \sim \sigma^{-1}$ leads to the desired outcome. Thus, combined with Theorem~A, ERM achieves the minimax rate $(s_N^*(\eta))^2$ for the confidence established in Theorem A (up to the constants in the exponent).
\vskip0.3cm
The proof of  Theorem \ref{thm:lower-bound-general} requires several preliminary steps.
\vskip0.3cm
Let $\bX=(x_i)_{i=1}^N\in \Omega^N$ and consider the conditional probability measure $\nu_{f,\bX}$ defined above. Put ${\cal A}_f=\tilde f_N^{-1}({\cal B}(f,\eps_N))$ and let ${\cal A}_f |\bX=\{y\in\R^N:(y,\bX)\in\cA_f\}$ denote  the corresponding fiber of ${\cal A}_f$ (see Figure~\ref{fig:minimax_lower_bound}).

\begin{Lemma} \label{lemma:conditioned-measure}
For every $f \in \cF$,
$$
Pr\big(\big\{\bX=(x_i)_{i=1}^N : \nu_{f,\bX}({\cal A}_f|\bX) \geq 1 - \sqrt{\delta_N} \big\}\big) \geq 1-\sqrt{\delta_N}.
$$
\end{Lemma}
\proof
Fix $f \in \cF$ and let $\rho(\bX) = \nu_{f,\bX}({\cal A}_f|\bX)$. Then,
$$
1-\delta_N \leq \nu_{f,\bX} \otimes \mu^N({\cal A}_f)=\E\rho(X_1,...,X_N).
$$
Since $\|\rho\|_{L_\infty} \leq 1$ and $\E \rho(\bX) \geq 1-\delta_N$, by the Paley-Zygmund Theorem (see Chapter~3.3 in \cite{MR1666908}), $Pr(\rho(\bX) \geq x)\geq (\E \rho(\bX)-x)/(1-x) \geq 1-\delta_N/(1-x)$ for every $0<x<1$. The claim follows by selecting $x=1-\sqrt{\delta_N}$.
\endproof
Observe that for every $f \in \cF$ and $\bX=(x_1,...,x_N)$,  $\nu_{f,\bX}$ is  a gaussian measure on $\R^N$ with mean $P_\bX f =(f(x_i))_{i=1}^N$ and covariance matrix $\sigma^2 I_N$.
\begin{Lemma} \label{lemma:distance-vs-measure}
Let $t\mapsto\Phi(t)=\Pro(g \leq t)$ be the cumulative distribution function of a standard gaussian random variable on $\R$.
Let $u,v \in \R^N$ and consider the two gaussian measures $\nu_u\sim \cN(u,\sigma^2 I_N)$ and $\nu_v\sim\cN(v,\sigma^2 I_N)$. If $A \subset \R^N$ is measurable, then
\begin{equation*}
\nu_v(A) \geq 1-\Phi\big(\Phi^{-1}(1-\nu_u(A))+\|u-v\|_{\ell_2^N}/\sigma\big).
\end{equation*}

\end{Lemma}

The main component in the proof of Lemma \ref{lemma:distance-vs-measure} is a version of the gaussian shift theorem.
\begin{Theorem} \cite{MR1652329} \label{thm:KuL}
Let $\nu$ be the standard gaussian measure on $\R^N$ and consider $B\subset\R^N$ and  $w\in\R^N$. If $H_+=\{x\in\R^N : \inr{x,w} \geq b\}$ is a halfspace satisfying that $\nu(H_+)=\nu(B)$, then $\nu(w+B) \geq \nu(w+H_+)$.
\end{Theorem}

\noindent{\bf Proof of Lemma \ref{lemma:distance-vs-measure}.}
Let $\nu$ be the standard gaussian measure on $\R^N$. A straightforward change of variables shows that
$$
\nu_u(A)= \nu\big((A-u)/\sigma\big) \mbox{ and } \nu_v(A)=\nu\big((A-v)/\sigma\big).
$$
Let $B=(A-u)/\sigma$, $w=(u-v)/\sigma$ and set $\nu(B)=\alpha$. Using the notation of Theorem~\ref{thm:KuL}, the corresponding halfspace is
$$
H_+=\{x: \inr{x,w/\|w\|_{\ell_2^N}} \geq  \Phi^{-1}(1-\alpha)\},
$$
and therefore, if $w^\perp \subset \R^N$ is the subspace orthogonal to $w$,
$$
w+H_+=\{(\lambda+1)w+ w^\perp : \ \lambda \geq \Phi^{-1}(1-\alpha)/\|w\|_{\ell_2^N}\}.
$$
Clearly,
$$
\nu(w+H_+)=Pr\big(g \geq\Phi^{-1}(1-\alpha)+\|w\|_{\ell_2^N} \big),
$$
and the claim follows from Theorem \ref{thm:KuL} and the definition of $w$.
\endproof

\noindent \textbf{Proof of Theorem~\ref{thm:lower-bound-general}.}
Let $\tilde f_N$ be a procedure that performs with accuracy  $\eps_N \leq d^2_{\cF}(L_2)/4$ and a confidence parameter $\delta_N$. Shifting $\cF$ if needed, and since $\cF$ is star-shaped around one of its points, one may assume that $u=0\in \cF$ and consider $v \in \cF$ for which $4\eps_N \leq \|v\|_{L_2(\mu)}^2 \leq 8\eps_N$. By Chebyshev's inequality,
$Pr\big(\norm{P_\bX v}_{\ell_2^N}^2 \geq 4N\|v\|_{L_2(\mu)}^2\big) \leq 1/4$, and thus, for $\bX=(X_i)_{i=1}^N$
in a set of $\mu^N$-probability at least $3/4$, $\|P_\bX v\|_{\ell_2^N} \leq c_1\sqrt{N}\|v\|_{L_2(\mu)}.$

Let
$$
\cA_0=\tilde f_N^{-1}({\cal B}(0,\eps_N)) \ \ {\rm and} \ \ \cA_v=\tilde f_N^{-1}({\cal B}(v,\eps_N)),
$$
which, by the choice of $v$, are disjoint. Since $\tilde f_N$ performs with accuracy  $\eps_N$ and has a confidence parameter $\delta_N$, $\nu_{0,\bX}\otimes \mu^N(\cA_0) \geq 1-\delta_N$ and $\nu_{v,\bX} \otimes \mu^N(\cA_v) \geq 1-\delta_N$. Applying Lemma \ref{lemma:conditioned-measure}, with $\mu^N$-probability at least $1-2\sqrt{\delta_N}$,
\begin{equation} \label{eq:prob-est}
\nu_{0,\bX}(\cA_0|\bX) \geq 1-\sqrt{\delta_N}, \ \ {\rm and} \ \  \nu_{v,\bX}(\cA_v|\bX) \geq 1-\sqrt{\delta_N}.
\end{equation}
Let $\Omega_0$ be the set of samples $\bX=(X_i)_{i=1}^N \subset \Omega^N$ for which
$\|P_\bX v\|_{\ell_2^N} \leq c_1\sqrt{N}\|v\|_{L_2(\mu)}$
and \eqref{eq:prob-est} holds. Hence, $Pr(\Omega_0) \geq 3/4-2\sqrt{\delta_N}$, and by Lemma \ref{lemma:distance-vs-measure} applied to the set $\cA_0|\bX$,
$$
\nu_{v,\bX}(\cA_0|\bX) \geq  1-\Phi\left(\Phi^{-1}(\sqrt{\delta_N})+\|P_\bX v\|_{\ell_2^N}/\sigma\right)=(*).
$$
Observe that if $\delta_N < 1/4$ then $\Phi^{-1}(\sqrt{\delta_N})<0$ and $|\Phi^{-1}(\sqrt{\delta_N})| \sim \sqrt{\log(1/\delta_N)}$. Moreover, if
$\|P_\bX v\|_{\ell_2^N}\leq \sigma|\Phi^{-1}(\sqrt{\delta_N})|$
then $(*) > 1/2$.

Since $\bX \in \Omega_0$, $\|P_\bX v\|_{\ell_2^N} \leq c_1\sqrt{N}\|v\|_{L_2(\mu)}$; therefore, if
$$
\|v\|_{L_2(\mu)} \lesssim \sigma \sqrt{\frac{\log(1/\delta_N)}{N}},
$$
it follows that $\nu_{v,\bX}(\cA_0|\bX) > 1/2$. On the other hand,
$\cA_0|\bX$ and $\cA_v|\bX$ are disjoint and $\nu_{v,\bX}(\cA_v|\bX)\geq 1-\sqrt{\delta_N}$, which is impossible if $\delta_N < 1/4$.

Thus,
$$
\|v\|_{L_2(\mu)} \gtrsim \sigma \sqrt{\frac{\log(1/\delta_N)}{N}},
$$
and by the choice of $v$,
$$
8\eps_N \geq \|v\|_{L_2(\mu)}^2 \gtrsim \sigma^2 \frac{\log(1/\delta_N)}{N},
$$
as claimed.
\endproof

Next, let us turn to the proof of Theorem B, which is a straightforward application of the next observation:
\begin{Theorem} \label{thm:lower-bound-general-multiple-hypothesis}
There exists an absolute constant $c_0$  for which the following holds. Let $\cF$ and $Y^f$ as above, and assume that $\tilde f_N$ is a procedure that performs with accuracy $\eps_N=a_N^2$ and has a confidence parameter $\delta_N \leq 1/4$. For any $\theta\geq 4$ and $f\in\cF$, if $\Lambda$ is a $2a_N$-separated subset of $\cF \cap (f+\theta a_ND)$ then
$$
\log |\Lambda| \leq c_0 N \left(\frac{\theta a_N}{\sigma}\right)^2.
$$
\end{Theorem}

\proof
Observe that if $a_N \geq (1/2)d_{\cF}(L_2)$ then $|\Lambda|=1$ and Theorem \ref{thm:lower-bound-general-multiple-hypothesis} is trivially true. Hence, one may assume that $a_N <(1/2)d_{\cF}(L_2)$.

Let $a=a_N$, set $D(f,r)=\{h \in \cF : \|f-h\|_{L_2(\mu)} \leq r\}$ and put
$\Lambda$ to be a maximal $2a$-separated subset of $\cF \cap (f+\theta a D)$ with respect to the $L_2(\mu)$ norm. Thus, $\{D(f,a): f\in\Lambda\}$ is a family of disjoint subsets of $\cF \cap (f+\theta a D)$.

Recall that for any $\bX=(x_1,\ldots,x_N)\in \Omega^N$,  $\cA_f|\bX$ is the fiber of $\cA_f=\tilde f_N^{-1}(D(f,a))$. Since $\tilde f_N$ performs with accuracy $a^2$ and has a confidence parameter $\delta_N=1-\alpha$, it follows that for any $f \in \Lambda$,
\begin{equation*}
 \E_\bX \nu_{f,\bX}(\cA_f|\bX)=\nu_{f,\bX}\otimes\mu^N(\cA_f)\geq\alpha.
\end{equation*}
If $u \neq v$ in $\Lambda$ and $A \subset\R^N$ then by Lemma~\ref{lemma:distance-vs-measure},
\begin{equation*}
  \nu_{u,\bX}(A)\geq 1-\Phi\big(\Phi^{-1}(1-\nu_{v,\bX}(A))+\norm{P_\bX v-P_\bX u}_{\ell_2^N}/\sigma\big).
\end{equation*}
Fix $v_0 \in \Lambda$. Since $\{\cA_v|\bX,v\in\Lambda\}$ is a family of  disjoint sets,
\begin{align*}
  1&\geq \sum_{v\in\Lambda}\nu_{v_0,\bX}(\cA_v|\bX)\geq \sum_{v\in\Lambda}\Big(1-\Phi\big(\Phi^{-1}(1-\nu_{v,\bX}(\cA_v|\bX))+\norm{P_\bX v_0-P_\bX v}_{\ell_2^N}/\sigma\big)\Big)=\sum_{v\in\Lambda}\int_{z_\bX(v)}^\infty \varphi(x)dx,
\end{align*}
where $\varphi$ is a density function of a the standard gaussian $\cN(0,1)$ and
\begin{equation*}
  z_\bX(v)=\Phi^{-1}(1-\nu_{v,\bX}(\cA_v|\bX))+\norm{P_\bX v_0-P_\bX v}_{\ell_2^N}/\sigma.
\end{equation*}
Taking the expectation with respect to $\bX$,
\begin{equation}
  \label{eq:sum1}
  1\geq\sum_{v\in\Lambda}\E_\bX\int_{z_\bX(v)}^\infty \varphi(x)dx,
\end{equation}
and it remains to lower bound each expectation.

Recall that
$$
\E_\bX \nu_{v,\bX}\big((\cA_v|\bX)^c\big)\leq 1-\alpha\leq1/4,
$$
and by Chebyshev's inequality, $Pr\big(\nu_{v,\bX}(\cA_v|\bX)\geq 3/4\big)\leq 1/3$. Therefore, with $\mu^N$-probability at least $2/3$,
\begin{equation*}
  \Phi^{-1}\big(1-\nu_{v,\bX}(\cA_v|\bX)\big)=\Phi^{-1}\big(\nu_{v,\bX}\big((\cA_v|\bX)^c\big)\big)\leq \Phi^{-1}(3/4):=\beta.
\end{equation*}
Another application of Chebyshev's inequality shows that with $\mu^N$-probability at least $2/3$,
\begin{equation*}
\norm{P_\bX v_0-P_\bX v}_{\ell_2^N}\leq (3/2)\sqrt{N}\|v_0-v\|_{L_2(\mu)}\leq (3/2)\theta a \sqrt{N},
\end{equation*}
because $v \in D(v_0,\theta a)$. Therefore, with $\mu^N$-probability at least $1/3$,
$$
z_\bX(v)\leq \beta+(3/2)\sqrt{N}\theta a/\sigma
$$
and since $\beta+(3/2)\sqrt{N}\theta a/\sigma>0$,
\begin{equation*}
  \E_\bX\int_{z_\bX(v)}^\infty \varphi(x)dx\geq \frac{1}{3}\int_{\beta+(3/2)\sqrt{N}\theta a_N/\sigma}^\infty \varphi(x)dx\gtrsim \exp\Big(-\frac{c_2N\theta^2 a^2}{\sigma^2}\Big).
\end{equation*}Thus, by (\ref{eq:sum1}),
 $ 1\gtrsim |\Lambda|\exp\big(-c_3N\theta^2 a^2/\sigma^2\big)$, as claimed.
\endproof

We end this section with the proof of Theorem~C, which is presented for a random choice of $(X_1,...,X_N)$, though the proof for a fixed $(x_1,...,x_N)$ -- the so-called deterministic design, is almost identical. The idea is that if $\bX=(X_1,...,X_N)$ and $P_\bX f_1 = P_\bX f_2$, the two functions are indistinguishable on a sample $(X_i,Y_i)_{i=1}^N$ of $Y^{f_1}=f_1(X)+V$. Therefore, no procedure can perform with a better accuracy than the largest typical $L_2(\mu)$ diameter of the sets
$$
K(f,\bX)=\{h \in \cF : P_\bX h =P_\bX f\}.
$$

Fix $f \in \cF$ and for every sample $\bX$ let $\cD(f,\bX)$ be the $L_2(\mu)$-diameter of $K(f,\bX)$. Define an $\cF$-valued random variable $h^f$ as follows. Let $h^f_{1,\bX}$ and $h^f_{2,\bX}$ be almost $L_2(\mu)$-diametric points in $K(f,\bX)$, set $\delta$ to be a $\{0,1\}$-valued random variable with mean $1/2$, which is independent of $X$ and $V$, and put
\begin{equation}
  \label{eq:r.v.h}
  h^f=(1-\delta)h^f_{1,\bX}+\delta h^f_{2,\bX}.
\end{equation}

Note that for every realization of $\delta$, $h^f\in K(f,\bX)$ and $\cD(h^f,\bX)=\cD(f,\bX)$. Denote by $Pr_{X,V}$ (resp. $\E_{X,V}$) the probability distribution of (resp. expectation w.r.t.)  $(X_i,V_i)_{i=1}^N$. Let $I(A)$ be the indicator of the set $A$ and observe that for every realization of the random variable $\delta$,
\begin{align*}
  & \sup_{f\in \cF} Pr_{X,V} \left( \| \tilde{f}_N \left((X_i,f(X_i)+V_i)_{i=1}^N\right)-f\|_{L_2(\mu)} \geq \cD(f,\bX)/4 \right)
  \\
 \geq & \sup_{f \in \cF}  Pr_{X,V} \left(\| \tilde{f}_N \left((X_i,h^f(X_i)+V_i)_{i=1}^N\right) - h^f\|_{L_2(\mu)} \geq \cD(h^f,\bX)/4 \right)
 \\
 = & \sup_{f \in \cF}  Pr_{X,V} \left(\| \tilde{f}_N \left((X_i,h^f(X_i)+V_i)_{i=1}^N\right) - h^f\|_{L_2(\mu)} \geq \cD(f,\bX)/4 \right) =(*)
\end{align*}
because $h^f\in\cF$ and $\cD(f,\bX)=\cD(h^f,\bX)$.

For every $f \in \cF$ put
$$
A_1^f= \left\{\|\tilde{f}_N\left((X_i,h^f_{1,\bX}(X_i)+V_i)_{i=1}^N\right) - h^f_{1,\bX}\|_{L_2(\mu)} \geq \cD(f,\bX)/4\right\},
$$
and
$$
A_2^f= \left\{\|\tilde{f}_N\left((X_i,h_{2,\bX}(X_i)+V_i)_{i=1}^N\right) - h_{2,\bX}\|_{L_2(\mu)} \geq \cD(f,\bX)/4\right\}.
$$
Taking the expectation in ($*$) with respect to $\delta$,
\begin{align*}
\E_\delta (*) \geq
\sup_{f \in \cF}  \E_{X,V} \E_{\delta } I\left(\tilde{f}_N \left((X_i,h^f(X_i)+V_i)_{i=1}^N \right) - h^f\|_{L_2(\mu)} \geq \cD(f,\bX)/4 \right)
=  \sup_{f \in \cF} \E_{X,V} \frac{1}{2}(I(A_1^f)+I(A_2^f)).
\end{align*}
Note that for any sample $\bX$, $h^f_{1,\bX}(X_i)+V_i=h^f_{2,\bX}(X_i)+V_i$; therefore,
$$
\tilde{f}_N \left((X_i,h^f_{1,\bX}(X_i)+V_i)_{i=1}^N\right) = \tilde{f}_N \left((X_i,h^f_{2,\bX}(X_i)+V_i)_{i=1}^N\right) \equiv f_0.
$$
Since $h^f_{1,\bX}$ and $h^f_{2,\bX}$ are almost diametric in $K(f,\bX)$, either $\|h^f_{1,\bX}-f_0\|_{L_2(\mu)} \geq \cD(f,\bX)/4$ or $\|h^f_{2,\bX}-f_0\|_{L_2(\mu)} \geq \cD(f,\bX)/4$. Thus, $I(A^f_1)+I(A^f_2) \geq 1$ almost surely, and
$$
\sup_{f\in \cF} Pr_{X,V} \left( \| \tilde{f}_N \left((X_i,f(X_i)+V_i,)_{i=1}^N\right)-f\|_{L_2(\mu)} \geq \cD(f,\bX)/4 \right) \geq 1/2.
$$
\endproof

\noindent{\bf Remark.} It is straightforward to verify that if $\sigma=0$, then ERM satisfies $\hat{f} \in K(f^*,\bX)$ for every sample $\bX$. Therefore, a typical value of $\cD(f^*,\bX)$ is a lower bound on the minimax rate when considering only noise-free targets.

As an example, let $T \subset \R^d$ be a convex, centrally-symmetric set, put $\mu$ to be an isotropic, $L$-subgaussian measure on $\R^d$ and set $\cF$ to be the class of linear functionals indexed by $T$. Given a sample $\bX=(X_1,...,X_N)$, set $\Gamma_{\bX}=\sum_{i =1}^N \inr{X_i,\cdot}e_i$ and put $P_\bX t =\Gamma_{\bX} t$. Therefore,
$$
K(v_0,\bX)=\{v \in T : \Gamma_{\bX} v = \Gamma_{\bX} v_0\} \subset 2T \cap {\rm ker}(\Gamma_{\bX}).
$$

Let $d_N=d_N(\rho)$ satisfy that with probability at least $1-\rho$, $\cD(0,\bX) \geq d_N$. Then, by Theorem~C, any procedure with a confidence parameter $\delta_N \leq 1/2+\rho$ cannot perform with a better accuracy than $d_N(\rho)/4$.

On the other hand, a straightforward application of Lemma \ref{lemma:above-r-N} shows that with probability at least $1-2\exp(-c_1NQ^2)$, $\cD(0,\bX) \lesssim r_N^*(Q)$. Therefore, if $d_N(T) \sim r_N^*(Q)$ for a suitable absolute constant $Q$, then with probability at least $1-2\exp(-c_1Q^2N)$,
$$
r_N^*(Q) \lesssim d_N(T) \leq \cD(0,\bX) \leq r_N^*(Q),
$$
and if $\sigma \lesssim r_N^*(Q)$, the error rate obtained in Theorem A is the minimax rate in the constant probability range.

\section{Examples} \label{sec:examples-TheoB}
In this section, we present two examples in which our results lead to sharp upper and lower minimax bounds, thus showing the optimality (in some minimax sense) of ERM.

\subsection{Learning in $\rho B_1^d$}
Let $\cF$ be the class of linear functionals $\inr{\cdot, t}$, indexed by $T=\rho B_1^d$, the unit ball in $\ell_1^d$ of radius $\rho$. Assume that $\mu$ is an isotropic, $L$-subgaussian measure on $\R^d$, that $Y \in L_{\psi_2}$ and that $\|Y-f^*(X)\|_{\psi_2} \leq \sigma$.

Since $\rho B_1^d$ is centrally symmetric, so is $\cF$, and $\cF-\cF = 2\cF$. Thus, the estimates in Theorem~A are based only on the behavior of the function $s \to \E\|G\|_{2\cF \cap s D}$. And, because the measure $\mu$ is isotropic, the canonical gaussian process is given by $t \to G_t = \sum_{i=1}^d g_i t_i$, where $g_1,\ldots,g_d$ are $d$ independent, standard gaussian variables. Moreover, for every $s>0$, the indexing set $2\cF \cap sD$ corresponds to $2\rho B_1^d \cap sB_2^d$. One may show (see, for example, \cite{MR2371614}) that for every $2\rho/\sqrt{d} \leq s$,
$$
\E \norm{G}_{2\rho B_1^d\cap s B_2^d} = \E \sup_{t \in 2\rho B_1^d \cap s B_2^d} \Big|\sum_{i=1}^d g_i t_i\Big| \sim \rho \sqrt{\log(ed \min\{s^2/\rho^2, 1\})},
$$
and if $s \leq 2\rho/\sqrt{d}$ then  $2\rho B_1^d \cap s B_2^d = sB_2^d$ and
$$
\E \norm{G}_{2\rho B_1^d\cap s B_2^d} = \E \sup_{t \in 2\rho B_1^d \cap s B_2^d} \Big|\sum_{i=1}^d g_i t_i\Big| \sim s\sqrt{d}.
$$
Setting $\eta = c_0/(L\sigma)$ and $Q=c_1/L^2$, it is straightforward to verify that
\begin{equation*}
(s_N^*(\eta))^2 \sim_L
\begin{cases}
\rho \sigma \sqrt{\frac{\log d}{N}} & \mbox{ if } \rho^2 N \leq \sigma^2 \log d, \hspace{1.5cm}  (c)\\
\rho \sigma \sqrt{\frac{1}{N}\log\Big(\frac{ed^2\sigma^2}{\rho^2N}\Big)} & \mbox{ if }  \sigma^2 \log d\leq \rho^2N \leq \sigma^2 d^2,
\\
\frac{\sigma^2 d}{N} & \mbox{ if } \rho^2 N \geq \sigma^2 d.
\end{cases}
\end{equation*}
Also,
\begin{equation*}
(r_N^*(Q))^2 \ \
\begin{cases}
\sim_L \frac{\rho^2}{N} \log\left(\frac{ed}{N}\right) & {\rm if} \ N \leq c_1d,
\\
\lesssim_L \frac{\rho^2}{d} & {\rm if} \ c_1 d \leq N \leq c_2d
\\
= 0 & {\rm if} \ N > c_2d,
\end{cases}
\end{equation*}
where $c_1$ and $c_2$ are constants that depend only on $L$.

When $N \sim d$, $(r_N^*(Q))^2$ decays rapidly from $(\rho^2/N)\log(ed/N)$ to $0$. Thus, when $c_1d \leq N \leq c_2d$ one only has an upper estimate on $(r_N^*(Q))^2$, and we will therefore only consider the cases $N \leq c_1d$ and $N \geq c_2d$.

\vskip0.3cm

Let us present the exact oracle inequalities satisfied by the ERM in $\rho B_1^d$ that follow from Theorem~A. First, assume that $N \leq c_1d$. If $\sigma \gtrsim r_N^*(Q)$ then $\sigma^2 d^2 \gtrsim N\rho^2$, and
$$
(s_N^*(\eta))^2 \sim_L
\begin{cases}
\rho \sigma \sqrt{\frac{\log d}{N}} & \mbox{ if } \rho^2 N \leq \sigma^2 \log d,\\
\rho \sigma \sqrt{\frac{1}{N}\log\Big(\frac{ed^2\sigma^2}{\rho^2N}\Big)} & \mbox{ if }  \sigma^2 \log d\leq \rho^2N.
\end{cases}
$$
Setting 
\begin{equation}\label{eq:proba_b1}
\delta_N = \left\{
\begin{array}{cc}
6\exp\left(-\frac{c_4\rho}{\sigma} \sqrt{N \log d}\right) &  \mbox{ if } \rho^2 N \leq \sigma^2 \log d,\\
6\exp\left(-\frac{c_4\rho}{\sigma} \sqrt{N \log\Big(\frac{ed^2\sigma^2}{\rho^2N}\Big)}\right) & \mbox{ if } \sigma^2 \log d\leq \rho^2 N,
\end{array}
\right.
\end{equation}
and applying Theorem A, it follows that if $\sigma \geq c_3 \rho \sqrt{\log(ed/N)/N}$,
then with probability at least $1-\delta_N$,
$$
R(\hat f) \leq \inf_{f \in \cF} R(f) + \frac{c_5\rho \sigma}{\sqrt{N}} \begin{cases}
 \sqrt{\log d} & \mbox{ if } \rho^2 N \leq \sigma^2 \log d,\\
\sqrt{\log\Big(\frac{ed^2\sigma^2}{\rho^2N}\Big)} & \mbox{ if }  \sigma^2 \log d\leq \rho^2N,
\end{cases}
$$
and if $\sigma \leq c_3\rho \sqrt{\log(ed/N)/N}$, then with probability at least $1-6\exp(-c_4N)$,
$$
R(\hat f) \leq \inf_{f \in \cF} R(f) + \frac{c_5\rho^2}{N} \log\left(\frac{ed}{N}\right),
$$
for constants $c_3,c_4,c_5$ that depend on $L$.

In a similar fashion, if $N \geq c_2d$ then $r_N^*=0$, and thus, if $\sigma \neq 0$,  $\sigma \geq r_N^*$. Therefore, the error rate of ERM is given by $s_N^*$. When $\sigma=0$ (the noise-free case) then $s_N^*(\eta)=r_N^*(Q)=0$ and with probability larger than $1-6\exp(-c_4N)$, $\hat f = f^*$, implying exact reconstruction.

\vskip0.3cm

Turning to the lower estimate, assume that the set of admissible targets contains every $Y^t=\inr{t,x}+W$, for $t \in \rho B_1^d$ and $W$ that is a centered gaussian random variable with variance $\sigma^2$ that is independent of $X$. It follows from Theorem A$^{\prime\prime}$ that if $\sigma \gtrsim r_N^*(Q)$, ERM is an optimal procedure in the following sense: it achieves the accuracy
$$
(s_N^*(c/\sigma))^2 \sim \rho \sigma \sqrt{(1/N) \log(ed^2 \sigma^2/(\rho^2 N))}
$$
if $\rho^2 N \geq \sigma^2 \log d$, and the accuracy 
$$
(s_N^*(c/\sigma))^2 \sim \rho \sigma \sqrt{(1/N) \log d}
$$
if $(\sigma^2/\log d)\lesssim\rho^2 N \leq \sigma^2 \log d$ (note that when $(\sigma^2/\log d)\gtrsim \rho^2 N$ then $\delta_N$ in \eqref{eq:proba_b1} is larger than $1$ and the probability estimate $1-\delta_N$ is negative).

For a minimax lower bound that holds with constant probability we shall apply Theorem~B. To that end, let us bound the covering numbers
$\log N(\rho B_1^d \cap 2rB_2^d, rB_2^d)$ from below. First note that
\begin{equation*}
N (\rho B_1^d\cap 2 r B_2^d, r B_2^d) = N (B_1^d\cap (2r/\rho)B_2^d, (r/\rho) B_2^d)
\end{equation*}
and it suffices to study the covering numbers $N(B_1^d\cap 2r B_2^d, r B_2^d)$ for various choices of $r$.

Fix $1/\sqrt{d} \leq 2r < 1$, and without loss of generality assume that $k=1/(2r)^2$ is an integer. For $I \subset \{1,...,d\}$, let $S^I$ be the Euclidean sphere supported on the coordinates $I$, and note that
$$
\bigcup_{|I|=k} 2rS^I \subset B_1^d \cap 2rB_2^d.
$$
It is a well known fact (see, e.g., \cite{MR2160043}) that there is a collection of subsets of $\{1,...,d\}$ of cardinality $k$, which will be denoted by ${\cal B}$, that is $k/8$ separated in the Hamming distance and for which $\log|{\cal B}| \geq c_1k\log(ed/k)$. Thus, the set $\Lambda = \{(2r)^2 \sum_{i \in I} e_i : I \in {\cal B}\}$ is an $r$-separated subset of $B_1^d \cap 2rB_2^d$ with respect to the $\ell_2^d$ norm, and for any $1/\sqrt{d}\leq 2r\leq1$,
$$
\log N(B_1^d \cap 2rB_2^d, r B_2^d) \geq c_4 \frac{\log(edr^2)}{r^2}.
$$
Moreover, one can prove (via Maurey's empirical method) that this estimate is sharp (see, e.g., \cite{MR732693}). Thus it follows that for any $\rho/\sqrt{d}\leq 2r\leq\rho$,
$$
\log N(\rho B_1^d \cap 2rB_2^d, r B_2^d) \sim  \frac{\rho^2}{r^2}\log\Big(\frac{edr^2}{\rho^2}\Big).
$$
If $2r \leq \rho/\sqrt{d}$ than $\rho B_1^d\cap 2r B_2^d = 2r B_2^d$ and by a volumetric estimate, $\log N(\rho B_1^d \cap 2rB_2^d, r B_2^d) \sim d$. If, on the other hand, $2\rho> 2 r \geq \rho$ then $\rho B_1^d\cap 2r B_2^d = \rho B_1^d$ and since $\log N(\rho B_1^d, r B_2^d) \sim \log(edr^2/\rho^2)\sim \log d$ (which is evident from the argument used above), then $\log N(\rho B_1^d \cap 2rB_2^d, r B_2^d) \sim \log d$. Finally, when $2r\geq 2 \rho$,  $\log N(\rho B_1^d \cap 2rB_2^d, r B_2^d)=0$.

Therefore,
 \begin{equation*}
(q_N^*(c_0/\sigma))^2 \sim_L
\begin{cases}
\rho^2 & \mbox{ if } \rho^2 N \leq \sigma^2 \log d, \hspace{1.5cm} (c^\prime)\\
\rho \sigma \sqrt{\frac{1}{N}\log\Big(\frac{ed^2\sigma^2}{\rho^2N}\Big)} & \mbox{ if }  \sigma^2 \log d\leq \rho^2N \leq \sigma^2 d^2,
\\
\frac{\sigma^2 d}{N} & \mbox{ if } \rho^2 N \geq \sigma^2 d.
\end{cases}
\end{equation*}

We conclude that when $\sigma \gtrsim r_N^*(Q)$ (and in particular, when $\sigma^2 d^2 \gtrsim N \rho^2$), and if $\rho^2 N \geq \sigma^2 \log d$, then  
$q_N^*(c_0/\sigma) \sim  s_N^*(\eta)$. This estimate also exhibits that Sudakov's inequality for the set $\rho B_1^d\cap2r B_2^d$ is sharp at the scale $\eps=r$ in the following sense: for every $0<r<\rho$,
\begin{equation*}
r \sqrt{\log N(\rho B_1^d\cap 2r B_2^d, r B_2^d)} \sim \E \norm{G}_{2\rho B_1^d\cap r B_2^d}.
\end{equation*}

Therefore, by Theorem~B, one has that if $\sigma \gtrsim r_N^*(Q)$ and if the set of admissible targets contains every $Y^t=\inr{X,t}+W$ as above, then the minimax rate in the constant confidence regime is $(s_N^*(\eta))^2$ and that ERM is optimal procedure when $\rho^2 N\geq \sigma^2 \log d$.

Note that when $\rho^2 N\leq \sigma^2 \log d$ the estimates on $(q_N^*(c_0/\sigma))^2$ and on $(s_N^*(\eta))^2$ do not coincide. And, it turns out that if one extends the set of admissible targets, ERM cannot perform with a better accuracy than $\sim (s_N^*(\eta))^2$ in this range. Indeed, consider the one dimensional case $d=1$ and a target $Y$ defined as follows: the marginal law of $Y$ given $X$ is
\begin{equation}\label{eq:other_stat_model}
Y = \left\{
\begin{array}{cc}
\sigma X & \mbox{ with probability  } 1/2 +\delta \\
-\sigma X & \mbox{ with probability } 1/2-\delta
\end{array}
\right.
\end{equation}
for $\delta$ that will be specified later, and $X$ that is distributed uniformly in $\{-1, 1\}$. The corresponding class of one-dimensional linear functionals is $\cF=\{f_t=tx \ : \ -\rho\leq t\leq \rho\}$.

It is straightforward to verify that for every $t\in[-\rho, \rho]$,
$$
R(t):=R(f_t) = (\sigma^2+t^2) - 4t\delta \sigma,
$$
and if $2\sigma \delta \geq \rho$ then the minimizer of $R(t)$ in $[-\rho,\rho]$ is $t=\rho$.

Next, let us identify the minimizer of the empirical risk $R_N(t) = N^{-1}\sum_{i=1}^N(Y_i-t X_i)^2$. Given the sample $(X_i,Y_i)_{i=1}^N$, let $J=\{i : Y_i = \sigma X_i\}$. Observe that for every $t \in [-\rho,\rho]$,
$$
\sum_{i=1}^N(Y_i-t X_i)^2=\sum_{i \in J} (\sigma-t)^2 X_i^2 +  \sum_{i \in J^c} (-\sigma-t)^2 X_i^2 = (\sigma -t)^2|J| + (\sigma+t)^2 |J^c|.
$$
Hence, if $|J^c| \geq |J|$ then the empirical minimizer satisfies $\hat{t} \leq 0$. By the choice of $Y$, the random variable $Z=\IND_{\{Y=-\sigma X\}}|X$ has mean $1/2-\delta$ and variance $\tau^2=1/4 - \delta^2$. Given $X_1,...,X_N$, let $Z_i=\IND_{\{Y_i=-\sigma X_i\}}|X_i$ and note that $|J^c| = \sum_{i=1}^N Z_i$.
It follows from the Berry-Esseen Theorem that if $\delta \sqrt{N}/\tau \leq c_1$ then with probability at least $1/4$, $|J^c| \geq N/2$. And to ensure that $\delta \sqrt{N}/\tau \leq c_1$ it suffices to select $\delta=c_2/\sqrt{N}$. All that remains is to estimate the excess risk of $\hat{t}$, which clearly satisfies
\begin{equation}\label{eq:lower_rate_model_shahar}
  R(\hat t) - R(t^*) \geq c_3 \sigma \rho\delta= \frac{c_4 \rho \sigma}{\sqrt{N}}.
\end{equation}
Thus, when $\sigma \delta \geq \rho$ (i.e., when $\rho^2 N \lesssim \sigma^2$), the best accuracy that ERM can achieve with constant probability is $\sim (s_N^*(\eta))^2$.

\vskip0.3cm

Finally, turning to the low noise regime ($\sigma\lesssim r_N^*(Q)$), one can show that the rate $(r_N^*(Q))^2$ is actually sharp. Recall that by Theorem~C it suffices to show that the Gelfand $N$-width of $\rho B_1^d$ satisfies $c_N(\rho B_1^d) \sim r_N^*$. By a result due to Garanaev and Gluskin \cite{MR759962}, when $d\geq N$ one has
$$
c_N(\rho B_1^d) \sim \rho \min \left\{1,\sqrt{\frac{\log(ed/N)}{N}}\right\},
$$
and $c_N(\rho B_1^d)=0$ when $d<N$. Therefore, $c_N(\rho B_1^d)\sim r_N^*(Q)$ when either $N\leq c_1 d$ or $N > c_2 d$. In particular, when $0 \leq \sigma \lesssim r_N^*(Q)$, the minimax rate is $(r_N^*(Q))^2$ and it is achieved by the ERM.

\subsection{Low-rank matrix inference via the max-norm}
In this section, the goal is to estimate the real-valued output $Y$ by a linear function of a low-rank (or approximately low rank) matrix. Since the rank is not a convex constraint, one may consider ``a convex relaxation'' given by the factorization-based norm
\begin{equation*}
\norm{A}_{max}=\min_{A=U V^\top}\norm{U}_{2\rightarrow \infty}\norm{V}_{2\rightarrow \infty}.
\end{equation*}
Let $\cB_{max}$ be the unit ball relative to that norm and set $\cF=\{f_A=\inr{\cdot,A}:A\in \cB_{max}\}$. Thus,
\begin{equation*}
  \hat A_N\in\argmin_{\norm{A}_{max}\leq 1}\frac{1}{N}\sum_{i=1}^N\big(Y_i-\inr{X_i,A}\big)^2.
\end{equation*}
A similar estimator has been studied in \cite{CaiZhou13} for $Y=\inr{A^*,X}+W$, a random vector $X$ that is selected uniformly from the canonical basis of $\R^{p \times q}$, a noise vector $W$ that is either gaussian or sub-exponential with independent coordinates, and matrices in $\cB_{max}$ with uniformly bounded entries.

Assume that $X$ is isotopic and $L$-subgaussian relative to the normalized Frobenius norm, and in particular, 
$$
\norm{\inr{X,A}}_{L_2}=(pq)^{-1/2}\norm{A}_F, \ \ \ \norm{\inr{X,A}}_{\psi_2}\leq L(pq)^{-1/2}\norm{A}_F.
$$
Let $A^*\in\argmin_{A\in\cB_{max}}\E (Y-\inr{A, X})^2$ be a minimizer of the risk in $\cB_{max}$ and set $ \sigma=\norm{Y-\inr{X,A^*}}_{\psi_2}$. Since $\cF$ is convex, the minimizer is unique and the conditions of Theorem A are satisfied.

To apply Theorem A, one has to estimate the fixed points $r_N^*(Q)$ and $s_N^*(\eta)$ for $Q$ that depends only on $L$ and $\eta \sim_L \sigma^{-1}$.

Let $B_F$ be the unit ball relative to the Frobenius norm. Since $X$ is isotropic, the relative $L_2$ unit ball is
$$
D=\{f_A:\E\inr{X,A}^2 \leq1\}=\{\inr{\cdot, A} :A\in\sqrt{pq}B_F\},
$$
and the corresponding gaussian process has a covariance structure given by
$$
\E G_{f_A} G_{f_B}=(pq)^{-1}\inr{A,B}=(pq)^{-1}{\rm Tr}(A^\top B).
$$
A simple application of Grothendieck's inequality (see, e.g., \cite{SS05}) shows that
\begin{equation*}
  {\rm conv}\big(\cX_{\pm}\big)\subset \cB_{max}\subset K_G {\rm conv}\big(\cX_{\pm}\big)
\end{equation*}
where $K_G$ is the Grothendieck constant and $\cX_{\pm}=\{uv^\top:u\in\{\pm1\}^p, v\in\{\pm1\}^q\}$; in particular, ${\rm diam}(\cB_{max},L_2) \sim 1$.

Let $\mathfrak{G}=(g_{ij})_{1\leq i\leq p:1\leq j\leq q}$ be a matrix with independent, centered gaussian entries with variance $(pq)^{-1}$. Thus, for every $s>0$,
\begin{align*}
  \E\norm{G}_{(\cF-\cF)\cap s D}=\E\sup_{A \in 2\cB_{max} \cap s\sqrt{pq}B_F}|\inr{\mathfrak{G},A}|\leq 2\E\sup_{A \in \cB_{max}}|\inr{\mathfrak{G},A}|
  \leq 2 K_G \E\sup_{A\in {\rm conv}(\cX_{\pm})}|\inr{\mathfrak{G},A}|.
\end{align*}
By standard properties of gaussian processes,
$$
\E\sup_{A\in {\rm conv}(\cX_{\pm})}|\inr{\mathfrak{G},A}| \lesssim \max_{A\in\cX_{\pm}}\frac{\norm{A}_F}{\sqrt{pq}}\sqrt{\log|\cX_{\pm}|} \lesssim \sqrt{p+q}.
$$
In the reverse direction, by Lemma~3.1 in \cite{CaiZhou13}, if
$$
\frac{1}{\min(p,q)} \lesssim s^2 \lesssim 1,
$$
then
\begin{equation} \label{eq:sudakov-matrix}
s \log^{1/2} N(\cB_{max}\cap s\sqrt{pq} B_F, s \sqrt{pq/2}B_F) \gtrsim \sqrt{p+q}.
\end{equation}
Hence, it follows from Sudakov's inequality that in that range of $s$,
$$
\E\norm{G}_{sD \cap (\cF-\cF)} \sim \sqrt{p+q},
$$
and
$$
(s_N^*(\eta))^2 \sim \sigma \sqrt{\frac{p+q}{N}}, \ \ \ (r_N^*(Q))^2 \sim \frac{p+q}{N},
$$
as long as both are smaller than $1$ and larger than $1/\min\{p,q\}$; that is, when $p+q\lesssim N \lesssim pq$, $p+q\lesssim \sigma^2 N$ and $\sigma^2 (p+q)\min(p,q)^2\gtrsim N$.

Applying Theorem A, if $\sigma \gtrsim_{Q,L} \sqrt{(p+q)/N}$ then with probability at least $1-2\exp(-c_1\sqrt{N(p+q)}/\sigma)$, ERM satisfies that
$$
\E(Y-\inr{\hat{A},X})^2 \leq \inf_{A \in \cB_{max}}\E(Y-\inr{A,X})^2 + c_2(Q,L)\sigma \sqrt{\frac{p+q}{N}},
$$
and if $\sigma \lesssim_{Q,L} \sqrt{(p+q)/N}$, then with probability at least $1-2\exp(-c_1N)$,
$$
\E(Y-\inr{\hat{A},X})^2 \leq \inf_{A \in \cB_{max}}\E(Y-\inr{A,X})^2 + c_2(Q,L)\frac{p+q}{N}.
$$

To see that the estimate is sharp in the minimax sense when $\sigma \gtrsim \sqrt{(p+q)/N}$ (and as long as $s_N^*, r_N^* \lesssim 1$, i.e., $\sigma \lesssim \sqrt{N/(p+q)}$), observe that Theorem A$^{\prime\prime}$ implies that ERM achieves the minimax rate for the confidence parameter $\delta_N =\exp(-c_1\sqrt{N(p+q)}/\sigma)$. Moreover, by Theorem B and \eqref{eq:sudakov-matrix}, any procedure with confidence parameter $\delta_N \leq 1/4$ has accuracy $\eps_N \gtrsim \sigma \sqrt{\frac{p+q}{N}}$, matching the upper bound.

\section{Concluding remarks} \label{sec:conc-rem}
Subgaussian classes are the first family of unbounded classes one is likely to consider, and it turns out that just like bounded classes, the study of subgaussian learning problems may be carried out using a two-sided concentration argument. Unfortunately, this is as far as concentration goes: the substantial technical machinery needed for the proof of Theorem A is not true beyond the subgaussian framework, and the analysis of more `heavy-tailed' problems requires a totally different machinery (see \cite{MR3367000,LWCG,LM_sparsity,LugMen}). Moreover, in more heavy-tailed situations, ERM does not attain the optimal accuracy/confidence tradeoff.

\vskip0.4cm

The results presented in this article are sharp in many cases but not in every case. First, in the `high probability' range, Theorem A$^{\prime\prime}$ shows that when $\sigma \gtrsim r_N^*$ the performance of ERM is optimal in the minimax sense. However, if $\sigma \lesssim r_N^*$, the estimate we present happens to be sharp only for $\sigma=0$ (when the error rate is a typical value of $\cD^2(f^*,\bX)$), or for $\sigma \sim r_N^*$, when the error rate is $\sim (r_N^*)^2$. This gap is filled (almost completely) in \cite{MenLvG}.

In the constant probability regime the picture presented here is even less complete. For example, in `noisy situations' -- when $\sigma\gtrsim r_N^*$, the upper bound of $(s_N^*(c/\sigma))^2$ is sharp only if it happens to be equivalent to $q_N^*(c/\sigma)$. Unfortunately, this is not even true even for $\cF = \{ \inr{t,\cdot} : t \in B_p^d\}$, when $1+1/\log d < p <2$. Again, this gap was addressed in \cite{MenLvG} -- at least when considering a class of admissible targets of the form $Y^f=f(X)+W$. 
\vskip0.3cm
The case of linear functional in $\R^d$ is a good indication to what our estimates give in general: if $X$ is $L$-subgaussian then when considering targets of the form $Y^t=\inr{X,t}+W$ for a centered gaussian variable $W$ that is independent of $X$, and $t \in T$, ERM achieves the accuracy
\begin{equation*}
\max\{(s_N^*(c/\sigma))^2,(r_N^*(Q))^2\}
\end{equation*}
as long as $T$ is convex and centrally-symmetric. No procedure can outperform this rate, say with confidence at least $3/4$ provided that:
\begin{enumerate}
\item   $q_N^*\log^{1/2} N(T\cap 2q_N^* B_2^d,q_N^* B_2^d)\sim \E\norm{G}_{2T\cap q_N^* B_2^d}$ -- meaning that there is no gap in Sudakov's inequality at scale $\eps= q_N^*$.
\item $c_N(T) \sim r_N^*(T)$ -- meaning that $\sqrt{N}c_N(T\cap r_N^*B_2^d)\sim \E\norm{G}_{T\cap r_N^*B_2^d}$, and there is no gap in the Pajor-Tomczak-Jaegermann estimate on the Gelfand $N$-width of $T$ (see \cite{MR845980}).
\end{enumerate}
Let us mention once again that a complete characterization of the minimax rate in this case was recently established in \cite{MenLvG}, and the optimal procedure happens to be a minor modification of ERM: it is ERM performed in an appropriate net in $T$.

The parameter $s_N^*$ may be compared with the fixed points used in \cite{vandegeer90,MR1240719,vdg93,van_de_geer00,MR2166554}. In all those cases, the fixed points are associated with Dudley's entropy integral for the localized class, rather than with the localized gaussian process; as such, the resulting bounds are always weaker than ours. For example, the results in \cite{MR1240719} which deal with the same situation as Theorem~A$^{\prime\prime}$ show that if the noise level is large enough and there is no gap in both Sudakov's AND Dudley's inequalities at the correct level (given by the fixed point), ERM is a minimax procedure in expectation. Theorem~A$^{\prime\prime}$ clearly improves that result.

\vskip0.4cm

Finally, although the importance of convexity may have been obscured by the Bernstein condition, a uniform Bernstein condition implies that the class is convex, at least if a nontrivial error rate is to be expected.

Indeed, observe that if $\cF \subset L_2(\mu)$ is closed but not locally compact in $L_2(\mu)$ then the minimax rate of $Y^f=f(X)+W$ does not tend to $0$ as the sample size tends to infinity. This is an immediate outcome of Theorem B and the fact that there is some $r>0$ and $f \in \cF$ for which $f+rD$ contains an infinite set that is $r/4$ separated in $L_2(\mu)$. Thus, one may restrict oneself to classes that are locally compact, and, in which case, one has the following:
\vskip0.6cm
\begin{Theorem} \label{thm:E}
Let $\mu$ be a probability measure and let $X$ be distributed according to $\mu$. If $\cF$ is a locally compact subset of $L_2(\mu)$, the following are equivalent:
\begin{enumerate}
\item[i)] for any real valued random variable $Y \in L_2$, the minimum of the functional $f \to \E(Y-f(X))^2$ in $\cF$ is attained. And, if $f^*$ is such a minimizer, then for every $f\in\cF$,
  \begin{equation}
    \label{eq:One-sided-Bernstein}
    \E\big(f(X)-f^*(X)\big)^2\leq \E\big((Y-f(X))^2-(Y-f^*(X))^2\big).
  \end{equation}
\item[ii)] $\cF$ is nonempty and convex.
\end{enumerate}
\end{Theorem}

\proof
If $\cF$ is a nonempty, closed and convex subset of a Hilbert space, the metric projection $Y \to f^*$ exists and is unique. By its characterization, $\inr{f(X)-f^*(X),Y-f^*(X)}\leq 0$ for every $f \in \cF$, and
\begin{align*}
  \E\big((Y-f(X))^2-(Y-f^*(X))^2\big)
=\norm{f(X)-f^*(X)}_{L_2}^2+2\inr{f^*(X)-Y,f(X)-f^*(X)}
\geq \norm{f(X)-f^*(X)}_{L_2}^2.
\end{align*}

In the reverse direction, if $\cF$ is locally compact, the set-value metric projection onto $\cF$ exists, and since it is $1$-Bernstein for any $Y$, the metric projection is unique. Indeed, if $f^*_1,f_2^*\in\cF$ are minimizers then by the Bernstein condition,
\begin{equation*}
  \norm{f_1^*(X)-f_2^*(X)}_{L_2}^2\leq B\E\big((Y-f_2^*(X))^2-(Y-f_1^*(X))^2\big)=0.
\end{equation*}

Thus, any $Y \in L_2$ has a unique best approximation in $\cF$, making $\cF$ a locally compact Chebyshev set in a Hilbert space. By a result due to Vlasov \cite{Vla}, (see also \cite{Deu}, Chapter 12), $\cF$ is convex.
\endproof

\appendix
\section{Additional proofs} \label{sub:proof_of_lemma_ref_lem_complexity_regularity}
First note that the canonical gaussian process we are interested in is a restriction of the isonormal process on $L_2(\mu)$ to a subset (see Section~12 in \cite{MR1932358}). In particular, it inherits the linearity of the isonormal process -- a fact we shall use below.

\noindent{\bf Proof of Lemma~\ref{lem:complexity-regularity}.}
Fix $s_1 > s_2 >0$ and $f,h \in \cF$. Assume that $s_2 \leq \|f-h\|_{L_2(\mu)} \leq s_1$ and observe that since $\cF-\cF$ is star-shaped around $0$ and $0< s_2/\|f-h\|_{L_2(\mu)} <1$, it follows that
$$
u=s_2 \frac{f-h}{\|f-h\|_{L_2(\mu)}} \in s_2 D \cap (\cF - \cF).
$$
Therefore,
\begin{equation} \label{eq:in-lemma-star}
G_{f-h} = \frac{\|f-h\|_{L_2(\mu)}}{s_2} G_u \leq (s_1/s_2) \sup_{w \in s_2 D \cap (\cF-\cF)} G_w.
\end{equation}
Since \eqref{eq:in-lemma-star} clearly holds if $\|f-h\|_{L_2(\mu)} \leq s_2$,  by taking the supremum over all possible choices of $f-h \in s_1 D \cap (\cF - \cF)$, 
\begin{equation*}
\sup_{w \in s_1 D \cap (\cF-\cF)} G_w\leq (s_1/s_2) \sup_{w \in s_2 D \cap (\cF-\cF)} G_w,
\end{equation*}
which is equivalent to $\psi(s_1)/s_1\leq \psi(s_2)/s_2$; therefore, $\phi$ is non-increasing on $(0,+\infty)$.

The two other parts of the claim can be established using a similar argument and their proofs are omitted. 
\endproof

\vskip0.5cm

\noindent{\bf Proof of Lemma \ref{lemma:r-vs-s}.} First, assume that $\sigma \geq (c/Q)r_N^*(Q)$. Let $r < r_N^*(Q)$ and note that by Lemma \ref{lem:complexity-regularity},
$$
\E\|G\|_{rD \cap (\cF-\cF)} \geq Qr \sqrt{N} = \frac{Q\sigma}{rc} \cdot \frac{c}{\sigma} r^2 \sqrt{N}.
$$
Hence, if $(Q\sigma)/rc \geq 1$ then $\E\|G\|_{rD \cap (\cF-\cF)} \geq \frac{c}{\sigma} r^2 \sqrt{N}$, implying that $r \leq s_N^*(c/\sigma)$. But $(Q\sigma)/rc \geq 1$ is equivalent to $\sigma \geq (c/Q)r$, which holds for any $r<r_N^*(Q)$.

For the reverse direction, let $\sigma \leq (c/Q)r_N^*(Q)$ and set $r > r_N^*(Q)$. Thus, by Lemma \ref{lem:complexity-regularity},
$$
\E\|G\|_{rD \cap (\cF-\cF)} \leq Qr \sqrt{N} = \frac{Q}{r}r^2 \sqrt{N}.
$$
Hence, if $Q/r \leq c/\sigma$ then $r \geq s_N^*(c/\sigma)$. But $Q/r \leq c/\sigma$ if $\sigma \leq (c/Q)r$, which clearly holds.
\endproof


\begin{spacing}{0.9}
\begin{footnotesize}
\bibliographystyle{plain}
\bibliography{biblio}
\end{footnotesize}
\end{spacing}

\end{document}